\documentclass{article}

\usepackage{arxiv}

\usepackage[utf8]{inputenc} 
\usepackage[T1]{fontenc}    
\usepackage{hyperref}       
\usepackage{url}            
\usepackage{booktabs}       
\usepackage{amsfonts}       
\usepackage{nicefrac}       
\usepackage{microtype}      
\usepackage{lipsum}		
\usepackage{graphicx}
\usepackage{natbib}
\usepackage{doi}
\newcommand{\ye}{Y_{3,n},\eta_{3,n}}
\newcommand{\yef}{Y_{3,4},\eta_{3,4}}
\newcommand{\yefi}{Y_{3,5},\eta_{3,5}}

\newcommand{\Aff}{\mbox{Aff}}

\usepackage{overpic}
\newtheorem{theorem}{Theorem}[section]
\newtheorem{lemma}[theorem]{Lemma}

\newtheorem{corollary}[theorem]{Corollary}
\newenvironment{proof}{\paragraph{Proof:}}{\hfill$\square$}

\title{Periodic Points of Ward-Veech Surfaces}

\author{\hspace{1mm}Benjamin J. Wright \\
	Department of Mathematics\\
	University of Wisconsin-Madison\\
	Madison, WI 53706 \\
	\texttt{bwright8@wisc.edu} \\}

\renewcommand{\shorttitle}{Periodic Points}

\hypersetup{
pdftitle={A template for the arxiv style},
pdfsubject={q-bio.NC, q-bio.QM},
pdfauthor={David S.~Hippocampus, Elias D.~Striatum},
pdfkeywords={Veech Surfaces, Ward Surfaces, Periodic Points},
}

\begin{document}
\maketitle

\begin{abstract}
	For a non-arithmetic Veech surface, it is known that the set points having finite orbit under the Veech group, called the set of periodic points, is finite. However, few examples of these periodic point sets have been computed. In what follows, we compute new examples of periodic point sets for a family of non-hyperelliptic Veech surfaces with unbounded complexity and consider applications to classifying the holomorphic sections of certain surface bundles as well as applications to the study of infinitely generated Veech groups. \end{abstract}


\section{Introduction}
\noindent A \textit{translation surface} is a pair $(X,\omega)$ where $X$ is a closed Riemann surface and $\omega$ is a non-trivial Abelian differential on $X$. A (possibly orientation-reversing) diffeomorphism of $X$ is said to be \textit{affine} with respect to $\omega$ if, away from the zeros of $\omega$, $\phi$ acts affinely in the coordinates of $X$ defined by integrating $\omega$. The set of all such diffeomorphisms forms a group called the \textit{affine group}, denoted $\Aff(X,\omega)$. A natural question to ask about this structure is: Which points of $X$ have finite orbit under the action of $\Aff(X,\omega)$? Such points are called \textit{periodic points}.\\

\noindent Because the affine group is generically finite (see \cite{Moeller09}), the set of periodic points of $(X,\omega)$ is generically $X$, but under certain hypotheses, the set of periodic points is known to be finite (see \cite{Moeller04}, \cite{ASENS_2003_4_36_6_847_0}). A few examples of such finite periodic point sets have been computed. \\

\noindent To introduce these examples, three more definitions are required: First, a translation surface $(X,\omega)$ is called a \textit{Veech surface} if $SL(X,\omega)$, the group of linear parts, or ``derivatives," of orientation-preserving maps of $\Aff(X,\omega)$, has finite co-volume in $SL_2(\mathbb{R})$. Second, a translation surface $(X,\omega)$ is called \textit{arithmetic} - or an \textit{origami} - if there is a ramified covering map $p:X\mapsto Y$ of a genus $1$ translation surface $(Y,\eta)$, ramified only over the identity element of $Y$, such that $p^\ast \eta = \omega$. If $(X,\omega)$ does not have the property of being an origami, it is called \textit{non-arithmetic}. More generally, a surface $(X,\omega)$ is called \textit{primitive} if it does not arise as the pull-back of a ramified covering $p:X\mapsto Y$, where $(Y,\eta)$ is a translation surface. Third, a point $x\in X$ is called a \textit{Weierstra{\ss} point} if there exists a meromorphic function $f$ on $X$, having only one pole at $x$ of order at most the genus of $X$. One can now state a theorem of M\"{o}ller (see \cite{Moeller04}): 

\begin{theorem}[M\"{o}ller] Let $(X,\omega)$ be a non-arithmetic genus $2$ Veech surface. Then the set of periodic points of $(X,\omega)$ is precisely the union of the vanishing set of $\omega$ with the set of Weierstra{\ss} points of $X$.
\end{theorem}

\noindent There are infinitely many surfaces satisfying the hypothesis of the above theorem (see \cite{genus2veech}). A similar result has been shown for another infinite family (see \cite{apisa2020periodic}): 

\begin{theorem}[Apisa, Saavedra, and Zhang] Let $(X,\omega)$ be the translation surface generated by playing billiards in a regular polygon (i.e. an algebraic curve with equation $y^2 + x^n = 1$ and differential $dx/y$). Then the set of periodic points of $(X,\omega)$ coincides with the set of Weierstra{\ss} points of $X$. 

\end{theorem}

\noindent It is worth noting that the examples noted so far are also all \textit{hyperelliptic} surfaces, that is, the surfaces $X$ admit holomorphic involutions $\sigma$ that act affinely such that $X/\langle \sigma\rangle$ is a Riemann sphere. And actually, a result of Gutkin, Hubert, and Schmidt (see \cite{ASENS_2003_4_36_6_847_0}) shows that if $(X,\omega)$ is hyperelliptic and its Veech group is generated by finite-order elements, then all Weierstra{\ss} points are periodic. (Notice that their theorem assumes the finite-order elements preserve orientation, but we extend the result to the case that the affine group is generated by [possibly] orientation-reversing maps in the following section.)\\

\noindent One of our main results says that this hyperelliptic condition is necessary. Namely: 

\begin{theorem} Let $Y_{3,4}$ be the (non-hyperelliptic, primitive) genus $3$ Riemann surface defined by the equation $$y^8 = (u-2)(u+1)^2$$ having Abelian differential $$\eta_{3,4} = (u + 1)du/y^7$$ Then $(Y_{3,4},\eta_{3,4})$ has $4$ periodic points.\\
\end{theorem}

\noindent Later, we will see the affine group of this surface may be generated by (orientation-reversing) torsion, and the set of periodic points will be identified precisely, making this apparently the first example computed of a finite set of periodic points of a non-hyperelliptic primitive Veech surface. But in particular, note that because $X$ has genus $3$, it has at least $8$ Weierstra{\ss} points (see, for example: \cite{book:3970}). Hence, as a corollary: 

\begin{corollary} There exist primitive translation surfaces whose affine groups are generated by torsion and whose Weierstra{\ss} points are not all periodic.

\end{corollary}

\noindent The proof of the theorem relies on a decomposition of $(Y_{3,4},\eta_{3,4})$ into regular polygons first given by Ward. Actually, Hooper gave decompositions into semi-regular polygons for a family of surfaces called the Bouw-M\"{o}ller surfaces, but this study focuses only on a subset of those known as Ward-Veech surfaces, defined in following sections. We also show: 

\begin{theorem} Let $n \geq 5$ be odd. The Riemann surface defined by the equation $$y^{2n} = (u-2)(u+1)^2$$ with translation structure given by the differential $$(u+1)du/y^{2n-1}$$ has exactly one periodic point which is not a singularity. In particular, such a surface has either two or four periodic points.
\end{theorem} 

and \begin{theorem} Let $n \geq 6$ be even. The Riemann surface defined by the equation $$y^{2n} = (u-2)(u+1)^2$$ with translation structure given by the differential $$(u+1)du/y^{2n-1}$$ has exactly three periodic points which are not singularities. In particular, such a surface has either four or six periodic points.
\end{theorem} 

\noindent Before arriving at the proofs, let us recall the special properties of the translation structures of these surfaces. 

\section{Geometry of Translation Surfaces}
\subsection{Cylinder Decompositions, Parabolics, and Veech Surfaces}

\noindent It is necessary to present several more definitions before preceding to study the Ward-Veech surfaces. So, let $(X,\omega)$ be any translation surface. A \textit{saddle-connection} on $(X,\omega)$ is a path $\gamma:I\mapsto X$ such that $\gamma(0),\gamma(1) \in Z(\omega)$, the zero set of the differential, and otherwise $\gamma(t)\not\in Z(\omega)$, and also such that $\gamma$ is a geodesic of the metric $\omega^2$. A direction $\theta\in S^1$ is said to be a \textit{periodic direction} of $(X,\omega)$ if there exists a saddle-connection on $(X,\omega)$ parallel to $\theta$. \\

\noindent And again, let $(X,\omega)$ be any translation surface. A \textit{cylinder decomposition} of $(X,\omega)$ in the direction of $\theta$ is a decomposition of $(X,\omega)$ into maximal cylinders foliated by closed geodesics in the direction of $\theta$, isometrically embedded in $(X,\omega)$ away from their boundaries, where all pairs of cylinders in the collection have commensurate moduli. \\

\noindent The following two theorems of Veech will be of great aid to the present study of periodic points (see \cite{Veech89}):

\begin{theorem}[Veech] Let $\theta$ be a periodic direction of a Veech surface. Then, $(X,\omega)$ admits an orientation-preserving affine map $\phi$ with parabolic linear part whose fixed direction is $\theta$.

\end{theorem}

\begin{theorem}[Veech] Let $(X,\omega)$ be any translation surface that admits an orientation-preserving affine map $\phi$ with parabolic linear part whose fixed direction is $\theta$. Then, $(X,\omega)$ admits a cylinder decomposition in the direction of $\theta$, and $\phi$ acts as a power of a Dehn twist on each cylinder.

\end{theorem}

\noindent The existence of such cylinder decompositions allows one to conclude that most points on a translation surface (with parabolics in its Veech group) cannot be periodic points. Namely, there is the following lemma, called the \textit{rational height lemma}, which has appeared in various places (see \cite{aprhl20},\cite{apisa2020periodic}):

\begin{lemma}[Rational Height Lemma] Let $(X,\omega)$ be any translation surface admitting an orientation-preserving parabolic affine map $\phi$ with fixed direction $\theta$. Let $\zeta$ stand for the direction perpendicular to $\theta$. Let $x\in C\subset X$ be a point in a cylinder $C$ of the decomposition in the direction of $\theta$. If the distance from $x$ to the boundary of $C$ along the direction of $\zeta$ is not a rational multiple of the height of $C$ (the $\zeta$-component of the distance between points on opposite boundaries), then $x$ has infinite order under the action of $\Aff(X,\omega)$.

\end{lemma}

\begin{proof} By rotating $(X,\omega)$, it is enough to show the lemma for when the fixed direction of $\phi$ is horizontal. Now, if the vertical distance from $x$ to the boundary of $C$ is $0$, the statement is trivial. So, let $x$ be away from the boundary of the cylinder, and say $C$ has width $w$ and height $h$. Let $R$ be the rectangle $[0,w]\times(0,h)$. The universal cover of $C$ is $(-\infty,\infty)\times(0,h) = S$ and has $R$ as (the closure of) a fundamental domain. Let $(x_1,x_2)$ be the coordinates of a point in $R$ corresponding to $x$. The deck group is horizontal translations by integer multiples of $w$.\\

\noindent The action of $\phi$ on $C$ lifts to a map $U$ of $S$, and can be expressed in coordinates as $$ U= \left[\begin{array}{cc} 1 & qw/h \\ 0 & 1\end{array} \right]$$ for some $q\in\mathbb{Q}$. Now, suppose that $x$ has finite order under the action of $\phi$, so $\phi^n (x) = x$, for some $n\in\mathbb{Z}$. Lifting to $S$ and comparing horizontal coordinates, one sees  $$x_1 + (nqw/h) x_2 = x_1 +mw $$ for some $m\in\mathbb{Z}$. It follows that $$x_2 = mh/nq$$ That is, $x_2$ is a rational multiple of $h$.

\end{proof}

\noindent There are three more technical tools used in the proofs in the main theorems of our study:

\begin{lemma} Let $(X,\omega)$ be any translation surface admitting a parabolic affine map $\phi$ with fixed direction $\theta$. Let $x\in X$ have finite orbit under $\phi$. Then, there is a cyclic ordering of the points of the $\phi$-orbit of $x$, namely $x_{[0]},\dots,x_{[n-1]}$, such that the $\omega^2$-distance between $x_{[i-1]}$ and $x_{[i]}$ is constantly $w/n$, where $w$ is the width of the cylinder where $x$ lives, and this distance is minimal between distinct points of the $\phi$-orbit of $x$. The points $\{x_{[i]}\}$ are said to be \textit{evenly distributed}.

\end{lemma}

\begin{proof} There is some $n\in\mathbb{N}$ such that $\phi^n(x)$ realizes the minimal distance $d$ between $x$ and another distinct point of the $\phi$-orbit of $x$. Because $\phi$ acts as a translation on lines on $(X,\omega)$ in the direction of $\theta$, $d$ is the length between $\phi^{nk}(x)$ and $\phi^{nk + n}(x)$ for any $k$. Note that no members of the $\phi$-orbit of $x$ appear between any $\phi^{nk}(x)$ and $\phi^{nk + n}(x)$, or else, by applying $\phi^{-nk}$, one finds a point of the $\phi$-orbit of $x$ between $x$ and $\phi^n(x)$, contradicting the minimality of $d$. Now, for some minimal $K > 1$, $\phi^{nK}(x) = x$. One uses the minmiality of this $K$ to realize $d = w/K$. \\

\noindent So, set $x_{[i]} = \phi^{ni}(x)$. The above paragraph shows this is well-defined, exhausts the $\phi$-orbit of $x$, and realizes constant distance between its consecutive members.

\end{proof}

\noindent The above allows us to apply the lemmas below in the following sections:

\begin{lemma} Let a circle $S^1$ be partitioned into intervals $A$ and $B$ such that $A\cap B = \partial A = \partial B$. Fix $N > 0$ and such that $A$ has length at least $1/N$ times the circumference of $S^1$. If $Q$ is finite set of points that are evenly distributed on $S^1$, and $|Q| \geq N$, then $Q$ meets the interior of $A$. 

\end{lemma}

\begin{proof} There are points $q_0,q_n \in Q$ having respectively minimal and maximal arguments (taken in $[0,2\pi)$), $\theta_0$ and $\theta_n$. Number all the points of $Q$ by order of increasing argument (taken in $[0,2\pi)$). The counter-clockwise distance from $q_0$ to $q_1$ is the length $1/N$ times the circumference of $S^1$. But, this should be the same as the counter-clockwise distance from $q_n$ to $q_0$, which is at least the length of $A$, which has length strictly more than $1/N$ times the circumference of $S^1$. Therefore, some $q_i$ meets the interior of $A$.

\end{proof}

\noindent Lastly, there is a converse to the rational height lemma:

\begin{lemma} Let $(X,\omega)$ be a translation surface admitting a parabolic affine map $\phi$ with horizontal fixed direction and cylinder $C$ with height $h$ and width $w$ on which $\phi$ acts with $$D\phi = \left[ \begin{array}{cc} 1 & w/h\\ 0 & 1  \end{array} \right]$$ Let $p\in C$ have height $nh/d$ in the cylinder where $(n,d) = 1$. Then, the order of $p$ under $\langle\phi\rangle$ is $d$. 

\end{lemma}

\begin{proof} Let $x$ denote the horizontal coordinate of $p$ in the cylinder $C$, which is well-defined modulo $w$. Note that $\phi^d(p)$ has $x$-coordinate $x + hdn/(hd) w = x + nw$, which is the same as $x$ modulo $w$. Now, suppose that $\phi^m(p)$ had the same $x$-coordinate as $p$ modulo $w$ for some $m$. Then, $x = x + \frac{mnh}{hd}w + kw$ for some integer $k$. Subtracting $x$, dividing by $w$, and rearranging shows that $mn = -kd$. Because $n$ and $d$ are relatively prime, $d|m$.

\end{proof}

\subsection{The Action of a Matrix on a Translation Surface}

\noindent The basic elements of the theory of translation surfaces covered by this section are not necessary to introduce for the proof of the main theorems of our study, but they are important in proving a theorem of Gutkin, Hubert, and Schmidt and in studying the application of the main theorems. \\

\noindent The elementary idea is as follows: Given a translation surface $(X,\omega)$, one has a natural atlas of charts $$\mathfrak{U} = \{h_U(z) := \int_{x}^z \omega \, :\, U \mbox{ is a disk neighborhood of } x\} $$ and the domain of $h_U$ is $U$. Note chart transformations away from the singularity set are simply translations. Now, take an invertible, two-by-two matrix $A$ define a new atlas of charts $$\mathfrak{U}^A = \{h^A_U(z) := A\int_{x}^z \omega \, :\, U \mbox{ is a disk neighborhood of } x\} $$ where $A$ acts on a complex number as a vector in $\mathbb{R}^2$. Chart transitions in $\mathfrak{U}^A$ (away from singularities) are simply translations, and therefore, these charts define a complex structure $AX$. One verifies readily that the differential $A\omega$ is holomorphic with respect to the structure $AX$, and so, we have a new translation surface $(AX,A\omega) = A(X,\omega)$. Now, notice that $X$ and $AX$ may both be identified with a topological surface $S_g$ by maps that forget their underlying complex structures. In fact, the identification maps may be made identical because we did not change the topology of $X$ to form $AX$, only its charts. In particular, we have a homeomorphism $i_A:X\mapsto AX$, which one verifies is an affine diffeomorphism with constant derivative $A$ by computing the derivative in coordinates. \\

\noindent In what follows, we will want to relate the affine groups (and in particular, the involutions of the affine groups) of $(X,\omega)$ and $(AX,A\omega)$. Concretely, we will want to use that:

\begin{lemma} If $(X,\omega)$ is a hyperelliptic translation surface with affine involution $\sigma$, than $(AX,A\omega)$ is hyperelliptic with affine involution $i_A\circ\sigma\circ i_A^{-1}$. Further, the Weierstra{\ss} points of $X$ are identified with the Weierstra{\ss} points of $AX$ via the homeomorphism $i_A$.
\end{lemma}

\begin{proof} Because $\sigma$ acts on $(X,\omega)$ with derivative $-I$, $i_A\circ\sigma\circ i_A^{-1}$ acts on $(AX,A\omega)$ with conjugate derivative $A(-I)A^{-1} = -I$. In particular, $i_A\circ\sigma\circ i_A^{-1}$ acts holomorphically on $AX$, and the quotient $AX/\langle i_A\circ\sigma\circ i_A^{-1} \rangle$ is homeomorphic to $X/\langle \sigma\rangle$, a sphere via a map induced by $i_A$. \\

\noindent Because the fixed points of these involutions are the Weierstra{\ss} points of these hyperelliptic curves (see \cite{nla.cat-vn1156072}), we see clearly that these points are identified via $i_A$.


\end{proof}

\noindent 

\subsection{(Anti)conformal Maps and Weierstra{\ss} points} This section is concerned with proving the following theorem: \begin{theorem}[Gutkin, Hubert, Schmidt] Let $(X,\omega)$ be a hyperelliptic translation surface with affine involution $\sigma$ such that $\Aff(X,\omega)$ is generated by torsion. Then, the Weierstra{\ss} points of $X$ are periodic.

\end{theorem}

\noindent Note how this statement allows for the possibility that the torsion reverses orientation. Basically, the proof will be exactly the same as in \cite{ASENS_2003_4_36_6_847_0}, only that we need the additional observation:

\begin{lemma} Let $\rho$ be an orientation-reversing isometry of the translation surface $(X,\omega)$. $\rho$ preserves the set of Weierstra{\ss} points of $X$. 

\end{lemma}

\begin{proof} Let $W$ be the set of Weierstra{\ss} points of $X$ a curve of genus $g$, and let $x\in W$. So, there exists a meromorphic function $f:X\mapsto \mathbb{C}$ having its only pole at $x$ of order $n \leq g$. Because $\rho$ is anticonformal, the composition $\bar{f}\circ\rho$ is a meromorphic map having its only pole of order $n$ at $\rho^{-1}(x)$. Therefore, $\rho^{-1}(x)$ is a Weierstra{\ss} point of $X$. So, $\rho^{-1}(W)\subseteq W$. The same argument applied to $\rho^{-1}$ shows $\rho(W)\subseteq W$, which yields $\rho(W) = W$ overall.

\end{proof} 

\noindent Clearly, the same holds true for an orientation preserving isometry. A reminder on Weierstra{\ss} points on hyperelliptic curves: The Weierstra{\ss} points of a hyperelliptic curve are exactly the fixed points of the hyperelliptic involution, noted also a couple paragraphs above. These two points are both discussed in \cite{nla.cat-vn1156072}. We can now prove the theorem: 

\begin{proof}[Proof of theorem 2.8] Let $\rho\in \Aff(X,\omega)$ be torsion, and $D\rho$ be the derivative of $\rho$ as an affine map acting on $(X,\omega)$. From linear algebra, we know that there is an invertible matrix $A$ such that $AD\rho A^{-1}$ is a (possibly orientation-reversing) isometry of the plane with orthogonal columns. It is also the derivative of the (possibly orientation-reversing) affine symmetry $i_A \circ\rho\circ i_A^{-1}$ of the translation surface $A(X,\omega)$. Therefore, either $i_A \circ\rho\circ i_A^{-1}$ acting on $AX$ is conformal or anticonformal, but as noted, either way, $i_A \circ\rho\circ i_A^{-1}$ preserves the set of Weierstra{\ss} points of $A(X,\omega)$. By lemma 2.7, it follows that $\rho$ preserves the set of Weierstra{\ss} points of $(X,\omega)$. \\

\noindent In summary, the finite set of Weierstra{\ss} points of $(X,\omega)$ are preserved by the generating set of the affine group and therefore the entire affine group - meaning - they are periodic.

\end{proof}

\noindent Gutkin, Hubert, and Schmdit used these observations and methods very similar to methods below to compute the periodic points of the translation surface generated by billiards in the regular octagon (see \cite{ASENS_2003_4_36_6_847_0}). This computation was generalized by M\"{o}ller in \cite{Moeller04} and Apisa, Saavedra, and Zhang in \cite{apisa2020periodic} into different settings. Note that in \cite{apisa2020periodic}, the authors suggest applying some of these related techniques to the Ward surfaces, defined in the next section.

\section{The Ward-Veech Surfaces}
 \noindent \subsection{The Ward surfaces and their Veech groups} Let $Y_{3,n}$ denote the closed Riemann surface defined by $$y^{2n} = (u-2)(u+1)^2$$ on its affine part which has an Abelian differential $$\eta_{3,n} = (u + 1)du/y^{2n-1}$$ The translation surface $(Y_{3,n},\eta_{3,n})$ is called a \textit{Ward surface}. A Ward surface may also be called a \textit{Ward-Veech} surface because of the following theorem (see \cite{ward_1998}):

\begin{theorem}[Ward] A Ward surface is a Veech surface.

\end{theorem}

\noindent The Veech group of a Ward-Veech surface can be written explicitly as follows: Let $A,B,C,\mu, R$ be the following matrices: $$A = \left[\begin{array}{cc} -1 & -1 \\ 0 & 1 \end{array}\right] \quad B = \left[\begin{array}{cc} -1 & 2 \cos(\pi/n) \\ 0 & 1 \end{array}\right]\quad C =\left[\begin{array}{cc} 0 & -1 \\ -1 & 0 \end{array}\right]$$ $$\mu = \left[\begin{array}{cc} \csc(\pi/n) & -\cot(\pi/n) \\ 0 & 1 \end{array}\right] \quad R = \left[\begin{array}{cc} \cos(\pi/n) & -\sin(\pi/n) \\ \sin(\pi/n) & \cos(\pi/n) \end{array}\right]$$ and $I$ be the $2\times 2$ identity. Let $X^Y$ stand for $YXY^{-1}$. Then theorems of Ward (see \cite{ward_1998}) and Hooper (see \cite{hooper}) yield: 

\begin{theorem} $\Aff(\ye) \cong GL(Y_{3,n},\eta_{3,n}) = \langle A^\mu, B^\mu, C^\mu \rangle$ and $\Aff^+(\ye) \cong SL(Y_{3,n},\eta_{3,n}) = \langle (AB)^\mu, (BC)^\mu \rangle = \langle (AB)^\mu, R \rangle$, where the isomorphisms are given by differentiation.

\end{theorem} 

\noindent Above, note that $GL(X,\omega)$ stands for the group of derivatives of all affine maps of $(X,\omega)$, i.e., including the ones that reverse orientation, and $\Aff^+(X,\omega)$ denotes the orientation-preserving affine maps. Also, in what follows we will commonly use $\phi$ to denote the affine map with parabolic derivative $(AB)^\mu$ and $\psi$ to denote the affine map with elliptic derivative $R$. \\

\begin{figure}[t]
    \centering

    \hspace{-.4in}\begin{overpic}[scale = 1.4]{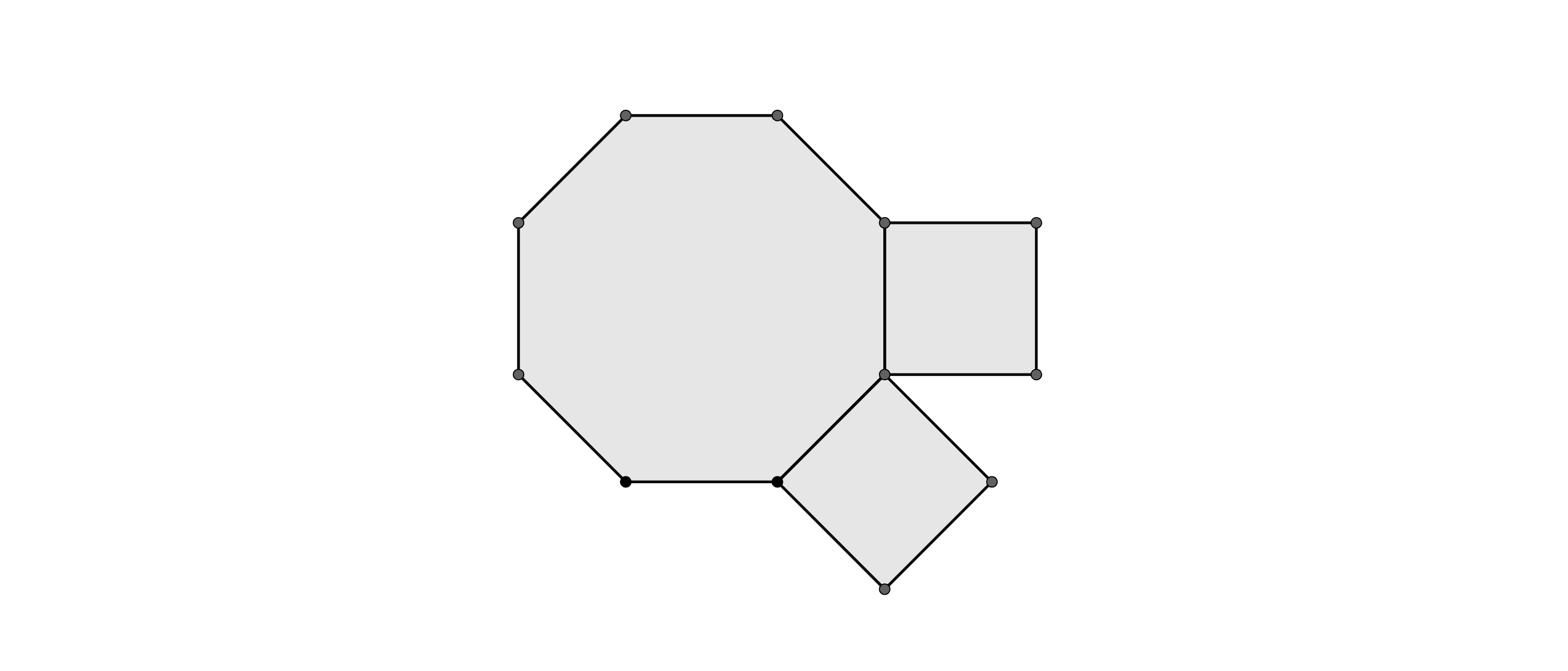}
    \put (30,23) {$a$}
    \put (67,23) {$a$}
    \put (34,33) {$b$}
    \put (61,7) {$b$}
    \put (44,37) {$c$}
    \put (64,17) {$c$}
    \put (54,33) {$d$}
    \put (50,7) {$d$}
    \put (60,30) {$e$}
    \put (44,10) {$e$}
    \put (34,13) {$f$}
    \put (62,15) {$f$}

    \end{overpic}
    \caption{The Ward-Veech surface $(\yef)$, whose pairs of edges with the same label are identified by translations}
    \label{fig:my_label}
\end{figure}
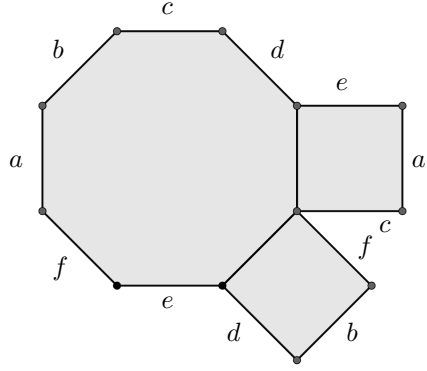

\noindent Later, we will also cite the following theorem of Hooper (see \cite{hooper}) about the translation structure of the Ward-Veech surfaces:

\begin{theorem} A Ward-Veech surface is a primitive translation surface.
\end{theorem}

\noindent \subsection{Ward surfaces via polygons} The Ward-Veech surfaces admit nice decompositions as gluings of regular polygons, originally due to Ward (see \cite{ward_1998}). Intuitively, the surface $(Y_{3,n},\eta_{3,n})$ may be seen as a gluing of three regular polygons all having the same side lengths. Two of the polygons are $n$-gons, and one is a $2n$-gon. The $n$-gons are rotated so that they may be glued by translation to any two sides of the $2n$-gon separated by an even number of edges. Furthermore, all edges of the $n$-gons are identified with edges in the $2n$-gon. (No polygon identifies a pair of its own edges.) Examples may be seen in Figures 1 and 4. We now recall the work of Hooper to make these constructions precise: Let $P_n(a,b)$ denote the $2n$-gon whose edge vectors $v_k$, $k = 0,\dots, 2n-1$ in order around the polygon are $$a(\cos (k\pi/n), \sin(k\pi/n))$$ if $k$ is even or $$b(\cos (k\pi/n), \sin(k\pi/n))$$ if $k$ is odd. Edges corresponding to even $k$ are called \textit{even}, and similarly, edges corresponding to odd $k$ are called \textit{odd}. Note that if exactly one of $a$ or $b$ is $0$, then $P_n(a,b)$ is a regular $n$-gon with only odd or even edges, respectively. Now, for $l = 0,1,2$ define $P(l)$ as: If $n$ is odd, $$P(l) = P_n(\sin((l+1)\pi/3),\sin(l\pi/3))$$ if $n$ is even and $l$ is even: $$P(l) = P_n(\sin(l\pi/3),\sin((l+1)\pi/3))$$ if $n$ is even and $l$ is odd, $$P(l) = P_n(\sin((l+1)\pi/3),\sin(l\pi/3))$$

\noindent It remains to describe the identifications between the edges of these three polygons. An even side of $P(1)$ is identified with the opposite side of $P(2)$, and an odd side of $P(1)$ is identified with the opposite side of $P(0)$, where \textit{opposite} means the identified sides have differing labels $k$ in their respective polygons. This construction is illustrated, for example, in Figure $1$. It is a theorem (see \cite{ward_1998}) that the translation surface arising via this construction is $(Y_{3,n},\eta_{3,n})$.

\subsection{Remarks on the horizontal cylinders of $(\ye)$.} Because the surfaces $(\ye)$ are Veech surfaces with horizontal periodic directions, they decompose into horizontal cylinders by theorem 2.1. Ward has shown (see \cite{ward_1998}) that the horizontal cylinders of $(\ye)$ all have identical moduli $$\alpha = (2\cos(\pi/n) + 1)/\sin(\pi/n)$$ One verifies that $((AB)^\mu)^{-1} = \left[ \begin{array}{cc} 1 & \alpha \\ 0 & 1 \end{array}\right]$ and therefore that lemma 2.6 applies to $\phi$ acting on any of the cylinders of this horizontal decomposition.
 
\section{Periodic Points of the Surface $(\yef)$} 

\noindent The surface $(\yef)$ can be illustrated as in the above Figure 1 as a gluing of three regular polygons. The form $\eta_{3,4}$ (arising via $dz$ in the figure) has a single zero of order $4$ at a vertex of the polygon, which can be seen by measuring the angle around this point. It follows that this surface has genus $3$.\\

\noindent Concretely, generators for the Veech group in this case are $$(AB)^\mu = \left[\begin{array}{cc} 1 & -2-\sqrt{2} \\ 0 & 1\end{array}\right]$$ and $$(BC)^\mu = \frac{\sqrt{2}}{2} \left[\begin{array}{cc} -1 & 1 \\ -1 & -1 \end{array}\right] $$ Noting that $((BC)^\mu)^4 = -I$, it is also possible to generate $SL(\yef)$ by $(AB)^\mu$ and the standard rotation of order $8$ $$R = ((BC)^\mu)^5 = -I(BC)^\mu = \frac{\sqrt{2}}{2} \left [\begin{array}{cc} 1 & -1 \\ 1 & 1 \end{array}\right] $$

\noindent The surface $(\yef)$ (being a Veech surface having horizontal and vertical saddle-connections, namely \textit{a,e} in Figure 1) admits horizontal and vertical cylinder decompositions, depicted in Figure 2 below. \\

\begin{figure}[b]
    \centering
    \hbox{\hspace{-.4in}\begin{overpic}[scale = 1.4]{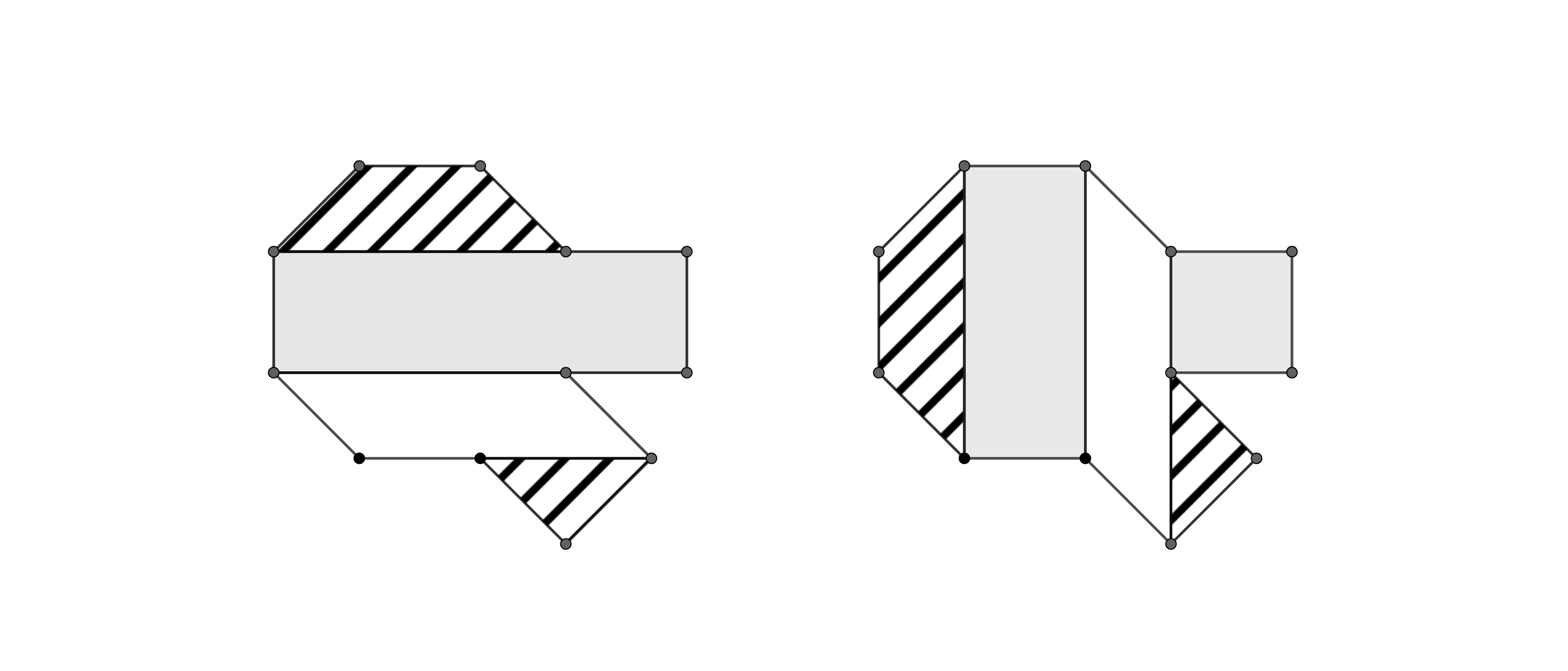}
    \put (34,30) {$H_1$}
    \put (45,22) {$H_2$}
    \put (41,16) {$H_3$}
    
    \put (56,15) {$V_1$}
    \put (64,11) {$V_2$}
    \put (72,30) {$V_3$}
    \end{overpic}}
    \caption{A decomposition of $(\yef) $ into $3$ horizontal and $3$ vertical cylinders - each with modulus $2 + \sqrt{2}$ }
    \label{fig:my_label2}
\end{figure}
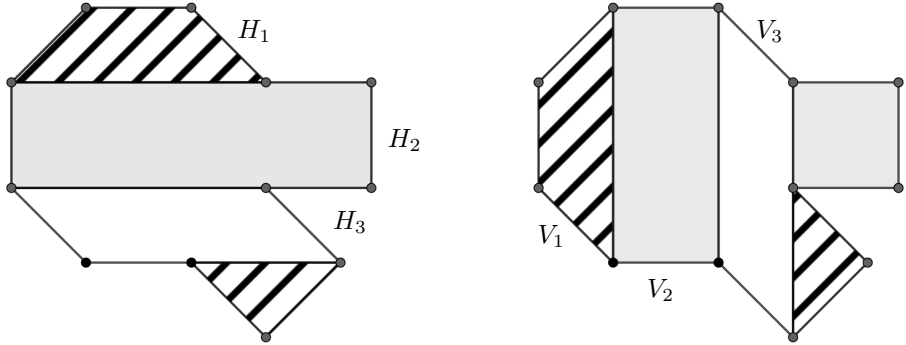

\noindent Theorem 2.1 guarantees a pair of parabolic affine maps, one for each of the cylinder decompositions, and we say a point in $y\in Y_{3,4}$ is a \textit{rational point} of the pair of cylinder decompositions if it has finite orbit under each of the corresponding pair affine maps. Note that by lemma 2.3, there are restrictions on the coordinates of a rational (and therefore, a periodic) point.\\

\noindent Again, let $\psi$ stand for the unique affine map of $(\yef)$ with derivative $R = -I(BC)^\mu$. The idea for the proof of the theorem is: If $y$ is a rational point, it is frequently the case that $\psi(y)$ is not, meaning that $y$ cannot be periodic. (It is worth noting that a point in a translation surface is periodic if and only if every point in its orbit under the affine group also is periodic - this will be used implicitly in what follows.) As another reminder, we let $\phi$ stand for the unique affine map with derivative $(AB)^\mu$. In the few cases that both $y$ and $\psi(y)$ are rational, it is often the case that there is an $n$ such that $\phi^n(y)$ belongs to the locus of $z$ where either $z$ or $\psi(z)$ is not rational, so $y$ is not periodic. It will remain to check that the finitely many points leftover from this process of elimination are periodic. One last note on notation: When we have fixed a chart in the plane of a translation surface $(X,\omega)$, we denote the $x$ and $y$ coordinates of a point $p\in X$ by $p_x$ and $p_y$, respectively. We are now ready for the proof.

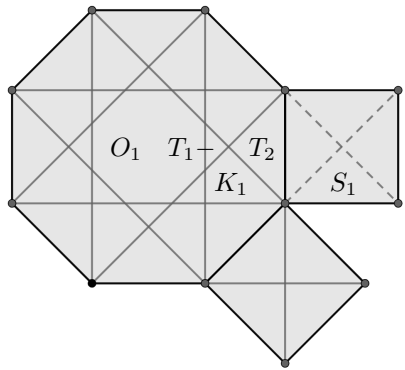
\begin{figure}[b]
    \centering
    \hbox{\hspace{0in} \begin{overpic}[scale = 1.5]{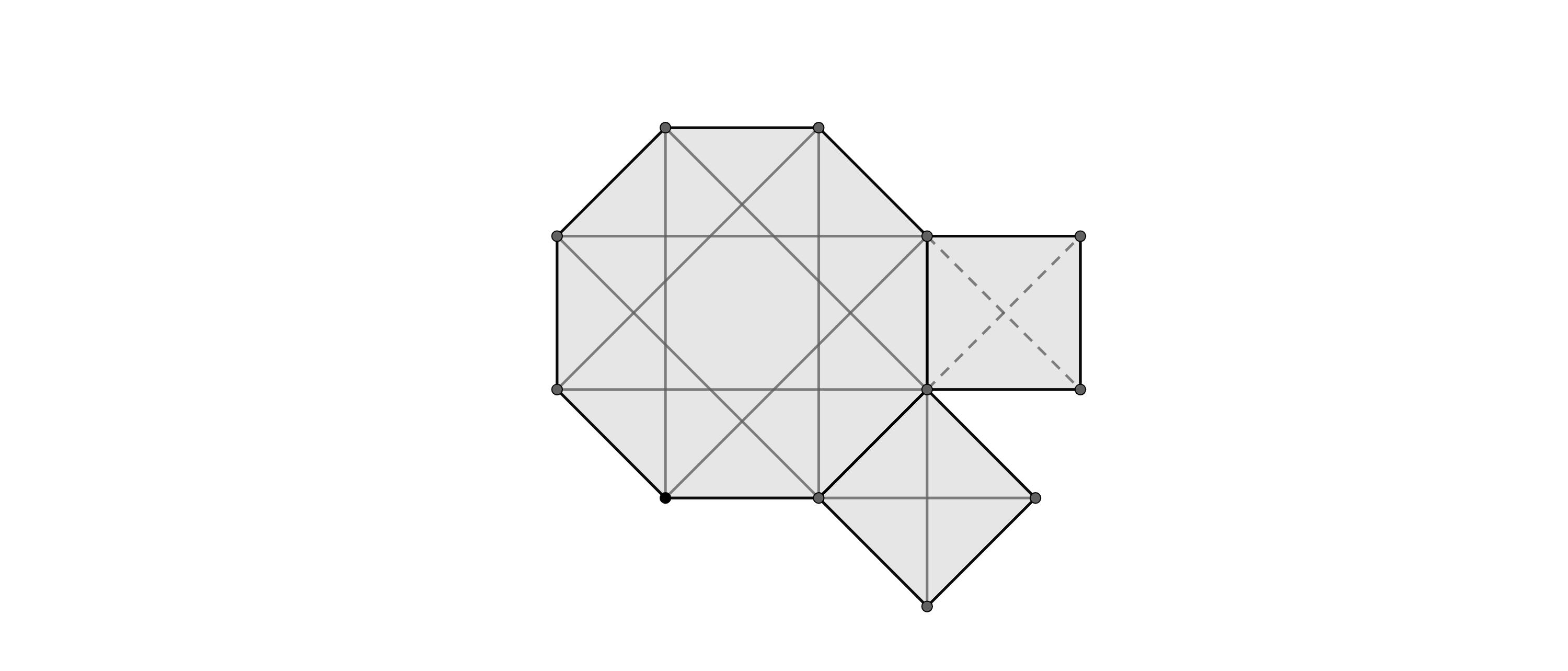}
    \put (44,22) {$O_1$}
    \put (49,22) {$T_1 - $}
    \put (56,22) {$T_2$}
    \put (53,19) {$K_1$}
    \put (63,19) {$S_1$}

    \end{overpic}}
    \caption{A decomposition of $(\yef)$ into polygons used to find periodic points.}
    \label{fig:my_label3}
\end{figure}

\begin{proof}[Proof of theorem 1.3] For convenience, change $(\yef)$ by scaling so that the length of the sides of the regular polygons forming $(\yef)$ (as in Figure 1) are $1$, and fix the center of the octagon at the origin by translating. For a point $p\in Y_{3,4}$, let $p_x,p_y$ be its $x$ and $y$ coordinates, respectively, in this coordinate system. Now, let $p\in Y_{3,4}$ such that $-1/2 \leq p_x,p_y \leq 1/2$, with the same holding true for $\psi(p)$. In other words, $p$ is on the interior octagon $O_1$ in Figure 3 (below). Suppose both $p$ and $\psi(p)$ are rational. Then, $p_x,p_y\in\mathbb{Q}$ by lemma 2.3 because $p\in H_2,V_2$. But also, $\psi(p)\in H_2,V_2$, so $$ \sqrt{2}/2(p_x-p_y),\sqrt{2}/2(p_x+p_y)\in\mathbb{Q} $$ by lemma 2.3. Adding and subtracting these expressions from each other gives $$\sqrt{2}p_x, -\sqrt{2}p_y\in\mathbb{Q}$$ which is not possible unless $p_x = p_y = 0$, i.e., $p$ is the origin. (It follows that the only possible periodic point in this case is the origin.)\\

\noindent Let $T_1$ be the triangle in Figure 3 and let $p\in T_1$ be such that $p,\psi(p)$ are rational. Hence, $p\in H_2,V_3$, but $\psi(p)\in H_2,V_2$. It follows that $$p_y,\psi(p)_x,\psi(p)_y,\sqrt{2}(p_x-1/2)\in\mathbb{Q}$$ by lemma 2.3. But then $\psi(p)_x -\psi(p)_y = -\sqrt{2}p_y \in\mathbb{Q}$. Hence, $p_y =0$. This means $\psi(p)_x = p_x/\sqrt{2}\in\mathbb{Q}$, in which case $\sqrt{2}(p_x-1/2)\not\in\mathbb{Q}$, a contradiction.\\

\noindent Now suppose $p,\psi(p)$ are rational, and $p\in K_1$, the kite (see Figure 3 below). Figure 3 shows that $p,\psi(p)\in H_2,V_3$. Hence, $p_y,\psi(p)_y\in \mathbb{Q}$ and $$\sqrt{2}(p_x - 1/2),\sqrt{2}(\psi(p)_x-1/2)\in\mathbb{Q}$$ by lemma 2.3. But $$p_x = \sqrt{2}\psi(p)_y - p_y$$ so the only way $\sqrt{2}(p_x - 1/2)\in \mathbb{Q}$ is if $p_y = -1/2$,  Combining this with the knowledge that $$\sqrt{2}(\psi(p)_x-1/2)\in\mathbb{Q}$$ gives that $p_x = q + \sqrt{2}/2$ for some rational $q$. Looking back to $\sqrt{2}(p_x-1/2)\in\mathbb{Q}$ shows that $p_x = 1/2 + \sqrt{2}/2 $. That is, $p$ is the vertex. \\

\noindent Now, let $p\in T_2$ be rational, and suppose $\psi(p)$ is also rational. Because $p\in H_2, V_3$ but $\psi(p)\in H_1,V_3$, it must be that $$p_y, \sqrt{2}(p_x -1/2), \sqrt{2}(\psi(p)_x - 1/2), \sqrt{2}(\psi(p)_y - 1/2)\in\mathbb{Q}$$ by lemma 2.3. It follows that $$ \sqrt{2}(\psi(p)_x + \psi(p)_y - 1) = (2p_x -\sqrt{2}) \in\mathbb{Q}$$ Therefore, $p_x = q + \sqrt{2}/2$ for a rational $q$. Because $\sqrt{2}(p_x - 1/2)\in\mathbb{Q}$, $q =1/2$. \\

\noindent In summary, up to the action of $\psi$, all points on the interior of the octagon, except the origin, have been ruled out as possible periodic points. \\

\noindent To finish narrowing down the set of periodic points, we translate $(\yef)$ vertically by +1/2, so that the bottom of the square $S_1$ lies on the $x$-axis. Let $p\in S_1$ be rational. Since $p\in H_2$, this means that $p_y\in\mathbb{Q}$ by lemma 2.3. Write $p_y$ as a reduced fraction $n/d$, where $n,d\in\mathbb{Z}$, $d\neq 0$. Also, assume that $p_y\neq 0, 1$, so it can be said that $d \geq 2$. Because $d,n$ are relatively prime, the orbit of $p$ by $\langle\phi\rangle$ has order $d$, by lemma 2.6. When $d\geq 3$, this orbit always intersects the interior of the octagon by lemmas 2.4 and 2.5. (To apply lemma 2.5, let $A$ be the intersection of the octagon with the closed leaf of the horizontal foliation passing through $p$ and $B$ be the intersection of the same leaf with the square containing $p$.) In the case $d = 2$, $n = 1$, and it is clear that each $p$ maps into the interior of the octagon because the image is contained in the component of $H_2\cap V_2$ that does not meet $S_1$. But note that the center of the square is mapped to the center of the octagon. \\

\noindent Therefore, up to the action of $\langle \phi,\psi\rangle$, all points except the singularity and the centers of the octagon and squares have been ruled out as possible periodic points. On the other hand, the named points are preserved  as a set by the group $\langle \phi,\psi\rangle$, which has index $2$ in $\Aff(\yef)$ - meaning these points are periodic.
\end{proof}

\subsection{Remarks on the Weierstra{\ss} points of $Y_{3,4}$} Knowing that $Y_{3,4}$ is given by $y^8 = (u-2)(u+1)^2$, one can directly compute that the Weierstra{\ss} weights of the $8$ points lying over $u = 5$ are all $2$. Also, the periodic points of $(\yef)$ are all Weierstra{\ss} points having weight $2$. But one does not need to know this to conclude corollary 1.3.1. 

\begin{proof}[Proof of Corollary 1.3.1] It is well-known that a curve of genus $g$ has at least $2g + 2$ Weierstra{\ss} points (see \cite{book:3970}). In this case $g = 3$, so $|P| < 2g+2 = 8$ proves the corollary.
\end{proof}

\noindent We also conclude:

\begin{corollary} The hyperelliptic hypothesis of theorem 2.8 is necessary. 

\end{corollary}

\begin{proof} The affine group of $(\yef)$ is generated by elliptic elements, but not all Weierstra{\ss} points are periodic. Therefore, by theorem 2.8, the surface is not a hyperelliptic translation surface. In other words: We find that the hyperelliptic condition is necessary for the theorem.

\end{proof}

\section{Periodic Points of the Surface $(\yefi)$}

\noindent Using similar methods, one can also compute the set of periodic points of the genus $4$ surface $(Y_{3,5},\eta_{3,5})$. In fact, we show the following: 

\begin{figure}[b]
    \centering
    \hbox{\hspace{.2in}\begin{overpic}[scale = 1]{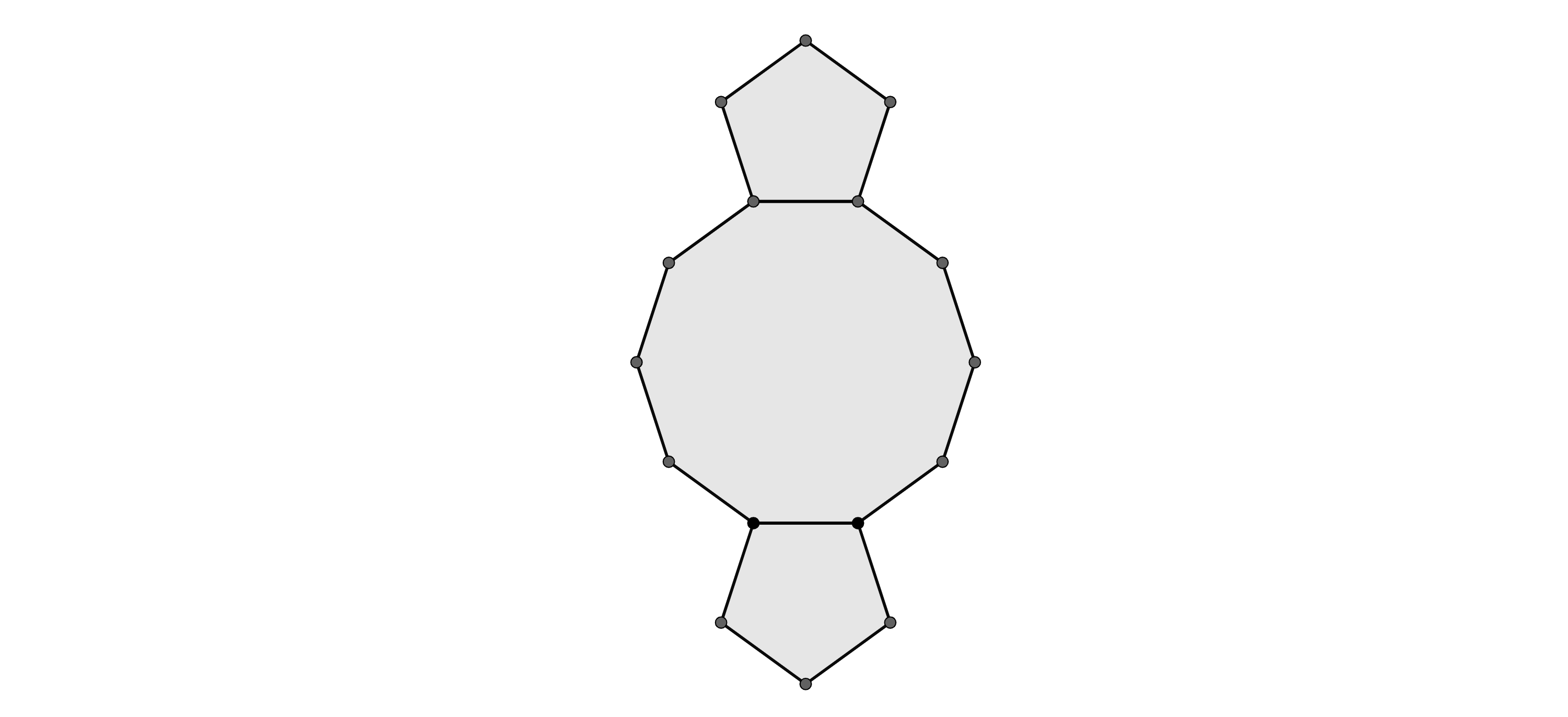}
    \put (47,42) {$a$}
    \put (58,13) {$a$}
    \put (54,42) {$b$}
    \put (43,13) {$b$}
    \put (57,36) {$c$}
    \put (39,25) {$c$}
    \put (44,36) {$d$}
    \put (62,25) {$d$}
    \put (45,10) {$e$}
    \put (62,19) {$e$}
    \put (56,10) {$f$}
    \put (39,19) {$f$}
    \put (47,3) {$g$}
    \put (57,32) {$g$}
    \put (54,2) {$h$}
    \put (44,32) {$h$}

    \end{overpic}}
    \caption{The Ward-Veech surface $(Y_{3,5},\eta_{3,5})$}
    
    \label{fig:my_label4}
\end{figure}
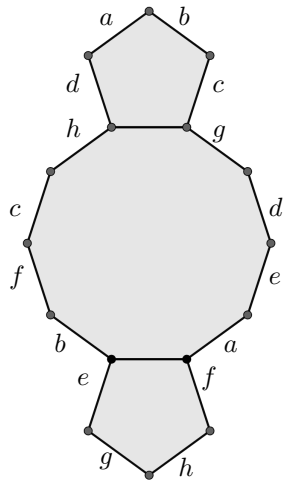

\begin{theorem} $(Y_{3,5},\eta_{3,5})$ has two periodic points, corresponding to $u = 2,\infty$ on $y^{10} = (u-2)(u+1)^2$. \\

\end{theorem}

\noindent Note this theorem is already implied by theorem 1.4. It is included for illustrative purposes.\\

\noindent Before beginning the proof, it will be helpful to highlight the horizontal cylinder decomposition of this surface. First, translate and scale $(\yefi)$ so that the center of the decagon is at the origin and so that the decagon is inscribed in the unit circle. $(Y_{3,5},\eta_{3,5})$ decomposes into four horizontal cylinders. The cylinders $C_2,C_4$, depicted in Figure 5, have height $\sin(\pi/5)$. The height of the other two cylinders is $\sin(2\pi/5) - \sin(\pi/5)$. Let $\phi$ be the parabolic affine map whose derivative is $(AB)^\mu$ (and preserves this decomposition). \\
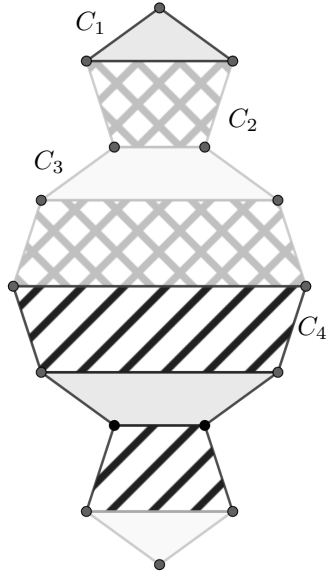
\begin{figure}[t]
    \centering
    \hbox{\hspace{-.1in}\begin{overpic}[scale = 1.2]{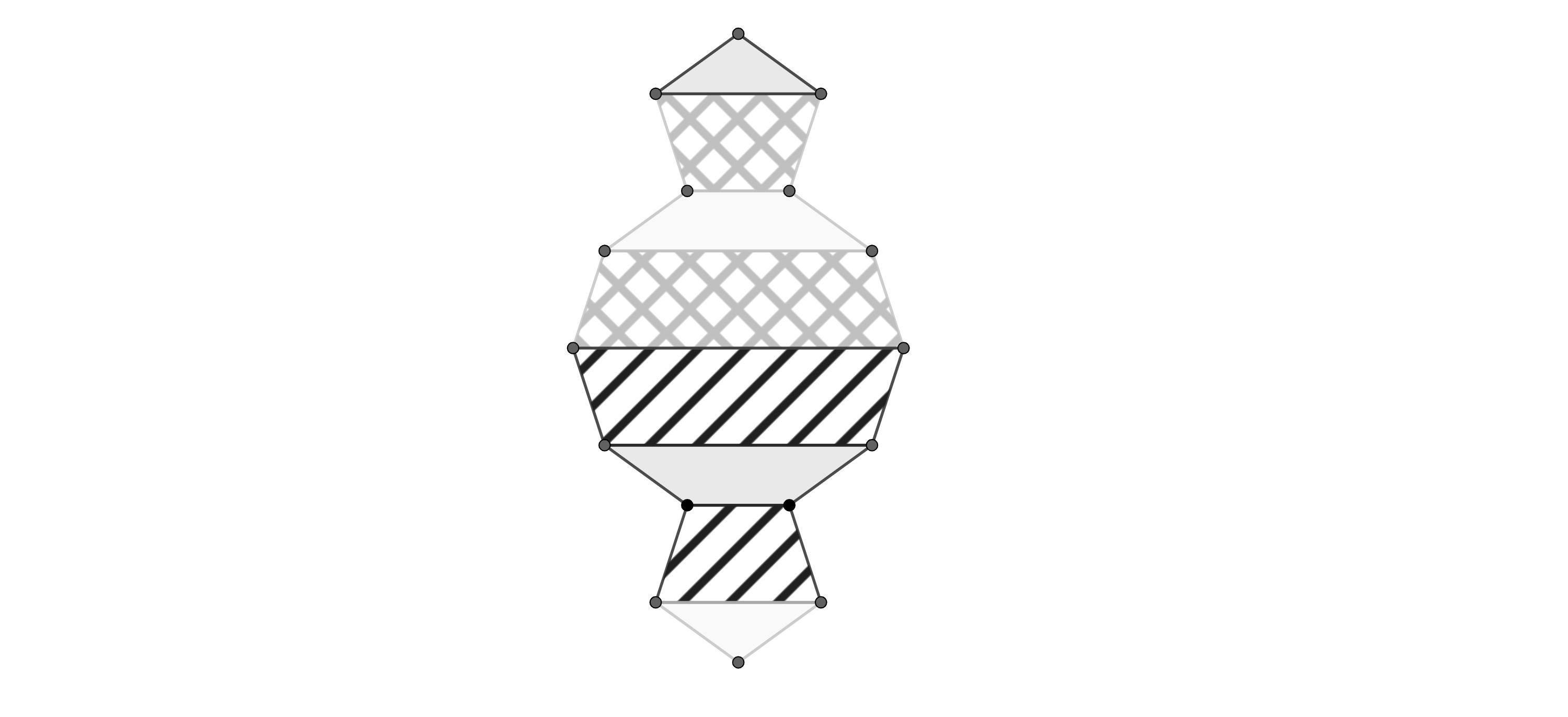}
    
    \put (41,42) {$C_1$}
    \put (52,35) {$C_2$}
    \put (38,32) {$C_3$}
    \put (57,20) {$C_4$}
    
    \end{overpic}}
    \caption{A horizontal cylinder decomposition of a genus $4$ Ward-Veech surface}
    \label{fig:my_label5}
\end{figure}

\noindent It will also be helpful to illustrate an order $10$ element of the affine group of $(Y_{3,5},\eta_{3,5})$, which will be called $\psi$ as usual. Apply the standard rotation of order $10$ to the polygon in Figure 4, centered at the center of the decagon. The two pentagons may be cut, translated, and re-glued according to the identifications along their boundary edges so that one again sees the polygon of Figure $4$. Recall from theorem 3.2 that $\Aff^+(\yefi) = \langle \phi, \psi \rangle$. We are now ready to begin the proof, which will be similar in spirit to the proof of theorem 1.3.\\
\begin{figure}[t]
    \centering
    \hbox{\hspace{-.55in}\begin{overpic}[scale = 2]{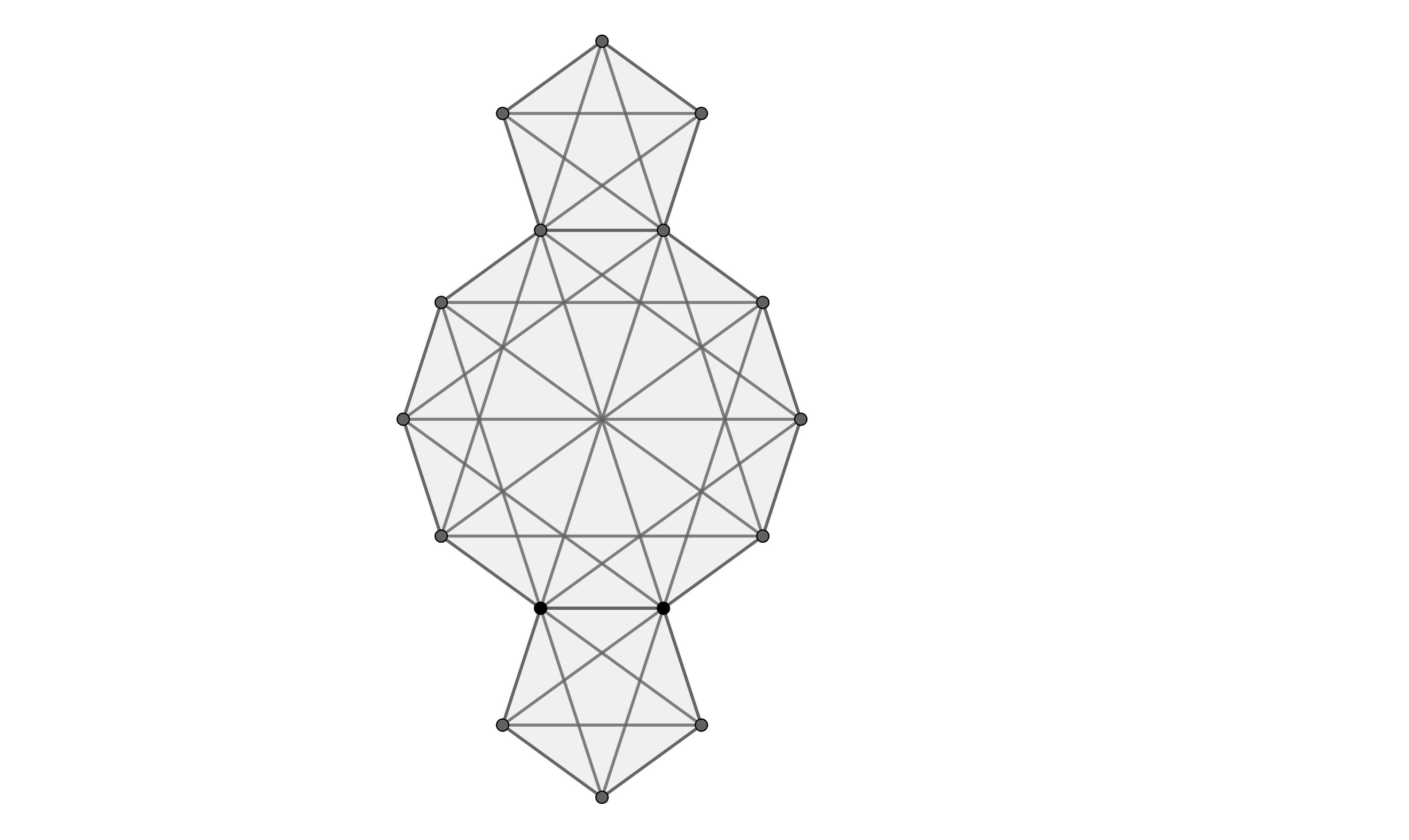}
    \put (47,31) {$T_1$} 
    \put (48.5,26.7) {$T_2 - $}
    \put (53,30.5) {$T_3$}
    \put (52.2, 28.2) {$T_4$}
    \put (42,41.5) {$T_5$}
    \end{overpic}}
    \caption{A helpful decomposition of $(Y_{3,5},\eta_{3,5})$ into polygons obtained by intersecting rotations of horizontal cylinders}
    \label{fig:my_label6}
\end{figure}
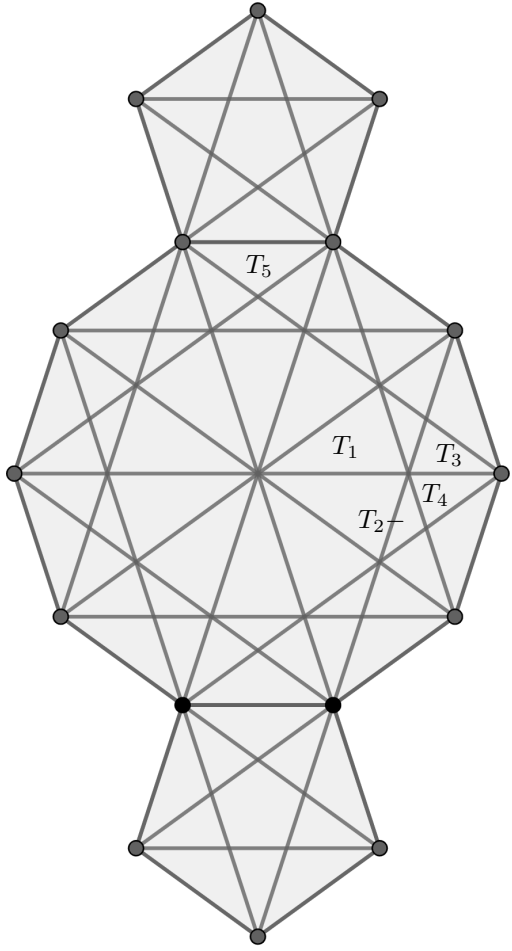

\begin{proof}[Proof of Theorem 5.1] Let $T_i$, $i = 1,2,3,4,5$, be the triangles depicted in Figure $6$. Assume $p\in Y_{3,5}$ is periodic. Now, let $p\in Y_{3,5} \in T_1,T_2,T_3$ or $T_4$. Then, $p$ has the property that $\psi^{-1}(p),p,\psi(p)\in C_2\cup C_4$. From $p\in C_2\cup C_4$, the rational height lemma shows that $p_y/\sin(\pi/5)\in\mathbb{Q}$. From $\psi^{-1}(p),\psi(p)\in C_2\cup C_4$, the rational height lemma also tells us that $$(\sin(\pi/5)p_x + \cos(\pi/5)p_y)/\sin(\pi/5),(-\sin(\pi/5)p_x + \cos(\pi/5)p_y)/\sin(\pi/5)\in \mathbb{Q}$$ Adding these  two expressions tells us that $$ 2\cos(\pi/5)p_y/\sin(\pi/5) \in \mathbb{Q}$$ But because $p_y/\sin(\pi/5)\in\mathbb{Q}$ and $\cos(\pi/5) = (1 + \sqrt{5})/4$ is irrational, it must be that $p_y = 0$. Notice this implies that $p_x\in \mathbb{Q}$ because $(\sin(\pi/5)p_x + \cos(\pi/5)p_y)/\sin(\pi/5)$ was already shown to be rational. Now, we must split into two sub-cases: \\

\noindent Case 1. First, suppose further that $p\in T_1$ or $T_2$. Then, $\psi^2(p)\in C_2$ also. Since $$\psi^2(p)_y = 2\sin(\pi/5)\cos(\pi/5) p_x$$ the rational height lemma says that $2\cos(\pi/5)p_x\in\mathbb{Q}$. But because $\cos(\pi/5)$ is irrational and $p_x\in\mathbb{Q}$, $p_x = 0$. In summary, the only periodic point in $T_1$ is the origin, and there are no periodic points in $T_2$.\\

\noindent Case 2. Now, the second sub-case: Suppose $p\in T_3$ or $T_4$. Then, $\psi^2(p)\in C_3$, and therefore the rational height lemma says $$(\psi^2(p)_y-\sin(\pi/5)) / (\sin(2\pi/5)-\sin(\pi/5))\in\mathbb{Q}$$ Using that $\psi^2(p)_y = 2\sin(\pi/5)\cos(\pi/5) p_x$ and simplifying with the double-angle formula for sine, we have  $$(2\cos(\pi/5) p_x -1)/(2\cos(\pi/5) -1)\in\mathbb{Q}$$ From the irrationality of $\cos(\pi/5)$ again, $p_x = 1$. Rephrased, the only periodic point in $T_3$ or $T_4$ is the vertex.\\

\noindent We have computed all possible periodic points in the union of $T_1,T_2,T_3,T_4$. Now, let $p\in T_5$. We also have that $\psi^{-1}(p),\psi(p)\in T_5$. The rational height lemma yields that 
$$\frac{p_y-\sin(\pi/5)}{\sin(2\pi/5)-\sin(\pi/5)},\frac{\sin(\pi/5)p_x + \cos(\pi/5)p_y-\sin(\pi/5)}{\sin(2\pi/5)-\sin(\pi/5)}, $$ $$\frac{-\sin(\pi/5)p_x +  \cos(\pi/5)p_y-\sin(\pi/5)}{\sin(2\pi/5)-\sin(\pi/5)} \in\mathbb{Q}$$

\noindent Adding the second two expressions and rearranging, one can write $$\cos(\pi/5)p_y - \sin(\pi/5) = t (\sin(2\pi/5) - \sin(\pi/5)) $$ for some $t,\in\mathbb{Q}$, and one can also write $p_y - \sin(\pi/5) = r ( \sin(2\pi/5) - \sin(\pi/5) ) $ for some $r\in\mathbb{Q}$. Subtracting these two expressions gives $$(\cos(\pi/5) - 1)p_y = u(\sin(2\pi/5) - \sin(\pi/5))$$ where $u = t-s\in\mathbb{Q}$. In particular, one has $$\frac{u(\sin(2\pi/5) - \sin(\pi/5))}{\cos(\pi/5)-1} = p_y = r(\sin(2\pi/5) - \sin(\pi/5)) + \sin(\pi/5)$$ or equivalently $$u(2\cos(\pi/5)-1) = r(2\cos(\pi/5)-1)(\cos(\pi/5)-1) + \cos(\pi/5) - 1$$ in which case one finds that $r = 1$ (and $u = -1/2$) using that $\cos(\pi/5) = (1 + \sqrt{5})/4$. Hence, $p_y = \sin(2\pi/5)$. This means the only periodic points in $T_5$ occur on the boundary of the decagon. \\

\noindent Now, up to the action of $\psi$, all points on the interior of the decagon except the origin have been ruled out as possible periodic points. It remains to show that no point on either of the pentagons, except the vertex, can be periodic. Note that for any point $p$ in either of the pentagons, $p$ has an image by a power of $\psi$ into either $C_2$ or $C_4$. Unless $p$ is the vertex, this image $\psi^k(p)$ can be taken to be away from the cylinder boundary. Now, all horizontal leaves of $C_2$ and $C_4$ have at least half of their lengths in the decagon, and therefore, by lemmas $2.4,2.5$, and 2.6, $\psi^k(p)$ has some image by a power of $\phi$ into the interior of the decagon. Notice that because the origin is on the cylinder boundary, this image cannot be the origin and therefore is not periodic.\\

\noindent It has been shown that the only possible periodic points of the surface $(\yefi)$ are the singularity and the origin. Because each of these points are fixed by $\phi,\psi$, they are fixed by the orientation-preserving affine group and are therefore periodic.

\end{proof}

 \section{Periodic Points of the Surface $(Y_{3,6},\eta_{3,6})$} 

 \noindent For the case $n = 6$, we have a similar result as in the case $n = 4$. Namely: 

\begin{theorem} The surface $(Y_{3,6},\eta_{3,6})$ has only $6$ periodic points: the singularities and the centers of the polygons (as depicted in Figure 7).
\end{theorem} 

\noindent Again, it is easy to see that these points are indeed left invariant (as a set) under the action of the affine group. The work will be in showing that no other points have finite orbits. The proof will be divided into cases based on the pattern of which cylinders of the horizontal decomposition a point in $(Y_{3,6},\eta_{3,6})$ visits under the action of the affine map whose derivative is the standard rotation of order 12 (calling this map $\psi$.) The cylinders of the decomposition are given labels according to Figure 8. 

\begin{figure}[t]
    \centering
    \includegraphics[scale = .5]{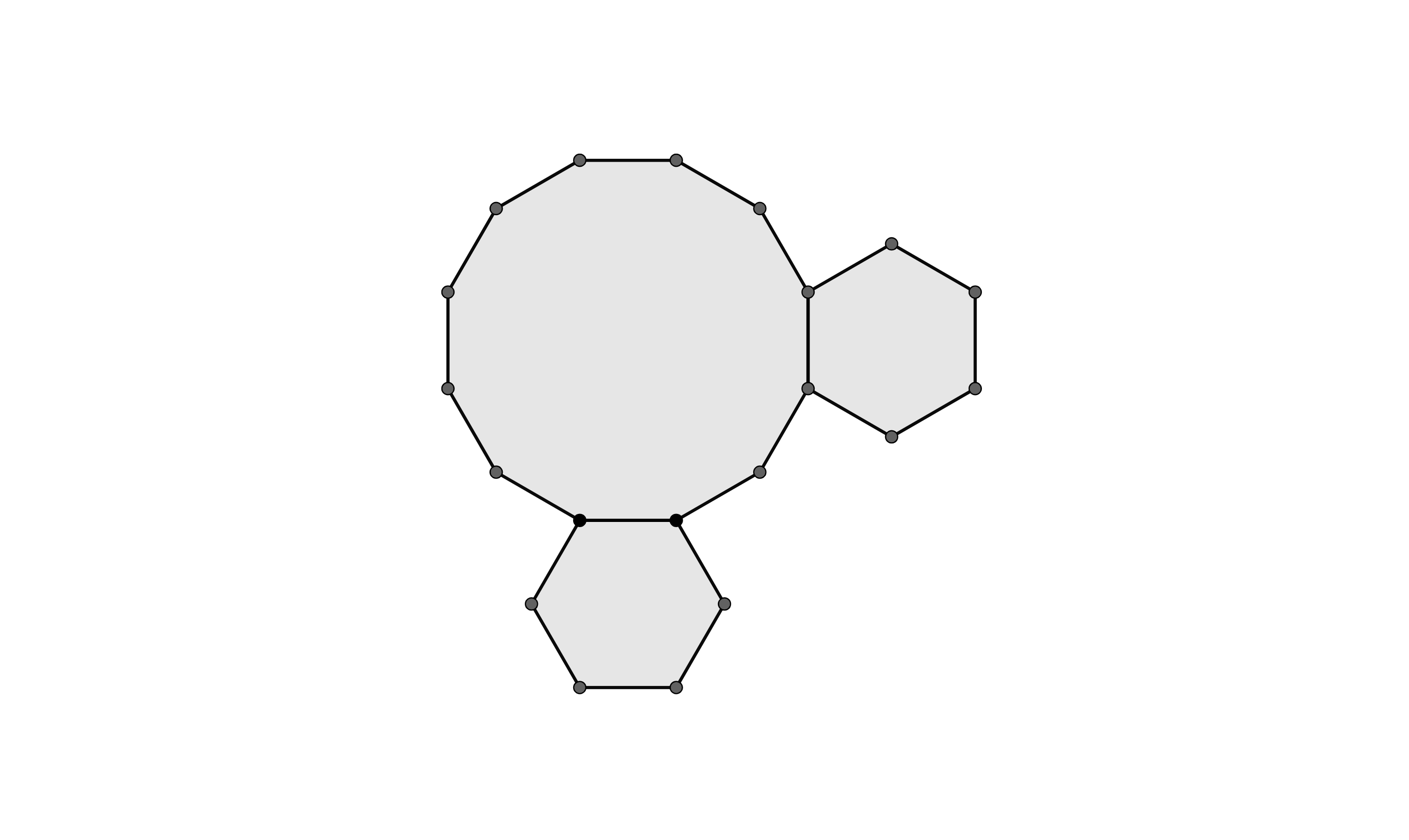}

    \caption{The surface $(Y_{3,6},\eta_{3,6})$. The edge identifications are explained in section $3.2$.}
    \label{fig:my_label7}
\end{figure}

\begin{figure}[t]
  
    \centering
    \hbox{\hspace{-1.2in}\begin{overpic}[scale = .8]{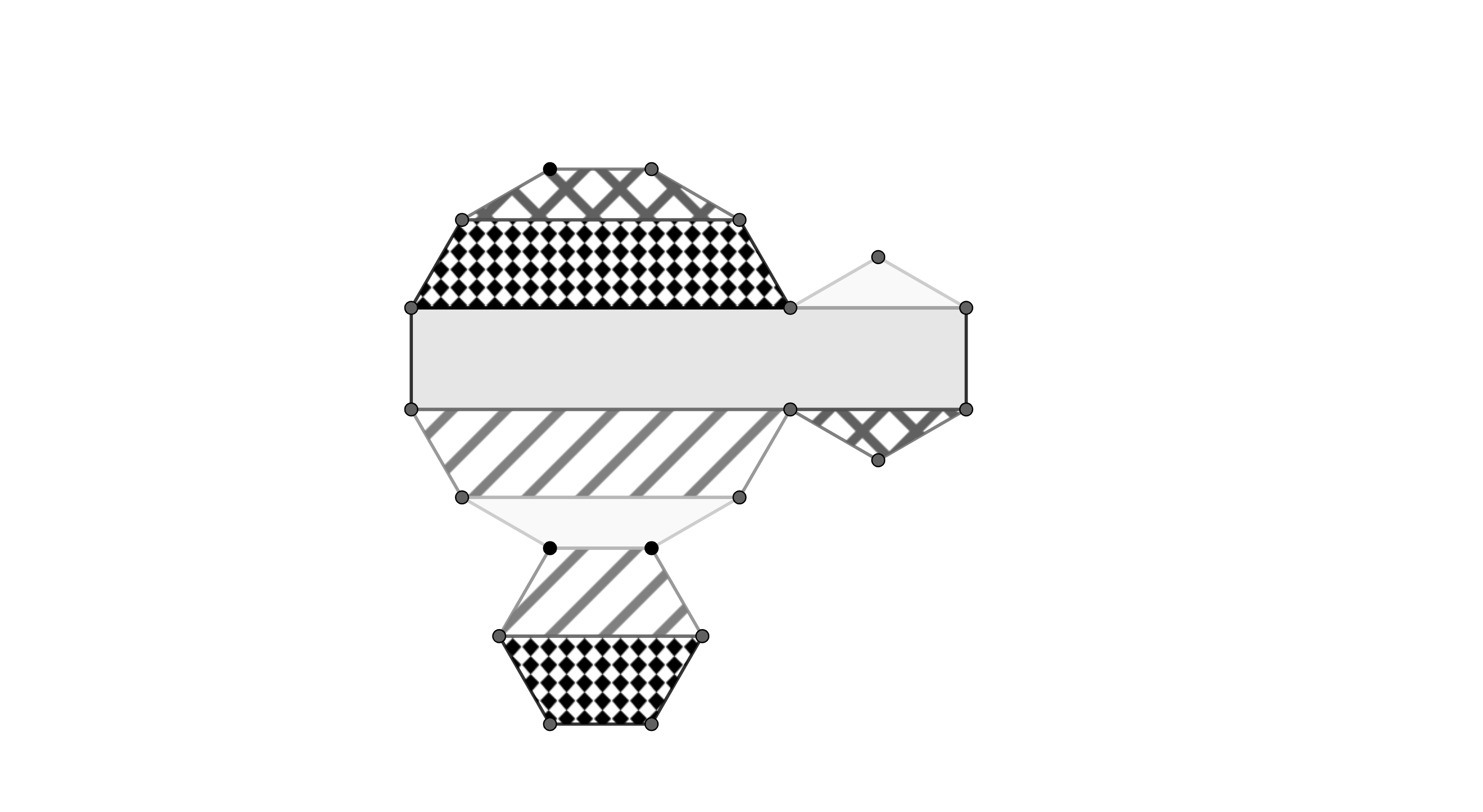}
    \put (49,42) {$C_3$}
    \put (25,37) {$C_2$}
    \put (66,30) {$C_1$}
    \put (52,23) {$C_4$}
    \put (32,18) {$C_5$}
    \end{overpic}}
    \caption{The horizontal cylinder decomposition of the surface $(Y_{3,6},\eta_{3,6})$}
    \label{fig:my_label8}
\end{figure}

\begin{proof}[Proof of Theorem 6.1] Scale and translate $(Y_{3,6},\eta_{3,6})$ so that the lengths of the polygonal sides are each $1$ and the center of the dodecagon is at the origin. Let $p$ be a point on the surface, and suppose it is periodic. \\

\noindent Case 1: Say $p$ has the property that $\psi^k(p) \in C_1$ for $k = -1,0,1$. The rational height lemma gives then that $$p_y, \sin(\pi/6)p_x + \cos(\pi/6)p_y, \sin(-\pi/6)p_x + \cos(-\pi/6)p_y \in \mathbb{Q}$$

\noindent Adding the second and third expressions gives $2\cos(\pi/6)p_y\in\mathbb{Q}$ which tells that $p_y = 0$ as $p_y\in\mathbb{Q}$. This means that $p_x\in \mathbb{Q}$ by substituting $p_y = 0$ back into the second expression. We break into sub-cases now. \\

\noindent Case 1a: Say $\psi^2(p)\in C_1$. Then, by the rational height lemma, $$2\sin(\pi/6)\cos(\pi/6)p_x = \cos(\pi/6)p_x\in\mathbb{Q}$$ which is impossible unless $p_x = 0$ because $p_x\in\mathbb{Q}$. Therefore, the only potential periodic point in this case is the origin.\\

\noindent Case 1b: Say $\psi^2(p)\in C_2$. By the rational height lemma, $$(2\sin(\pi/6)\cos(\pi/6)p_x - 1/2)/\sin(\pi/3) = p_x -1/\sqrt{3} \in \mathbb{Q}$$ This is impossible because $p_x\in\mathbb{Q}$.\\

\noindent Case 2: Say $p$ has the property that $p\in C_1$, $\psi^2(p)\in C_2$, and $\psi^{-2}(p)\in C_4$. So, the rational height lemma says $p_y\in\mathbb{Q}$ and $$(\sin(\pi/3)p_x + \cos(\pi/3)p_y - 1/2)/\sin(\pi/3), $$ $$ (-\sin(\pi/3)p_x + \cos(\pi/3)p_y + 1/2)/\sin(\pi/3) \in\mathbb{Q}$$ Adding these latter two expressions gives $p_y/\sin(\pi/3) \in\mathbb{Q}$, which implies $p_y = 0$ because $p_y\in\mathbb{Q}$. It follows that $$ p_x - 1/\sqrt{3} \in \mathbb{Q}$$ We now analyze sub-cases.\\

\noindent Case 2a: Say $\psi(p)\in C_1$. It follows that $$\sin(\pi/6)p_x = p_x/2\in \mathbb{Q}$$ This is impossible as $p_x - 1/\sqrt{3} \in \mathbb{Q}$, so there are no periodic points in this case.\\

\noindent Case 2b: Say $\psi(p) \in C_2$. It follows that $$(\sin(\pi/6)p_x - 1/2)/\sin(\pi/3) = (p_x - 1)/\sqrt{3} \in \mathbb{Q}$$ Because $p_x - 1/\sqrt{3} \in \mathbb{Q}$, we see that $p_x = 1 + 1/\sqrt{3}$. Suppose further that $\psi^3(p)\in C_2$. Then, the rational height lemma yields $$(p_x -1/2)/(\sqrt{3}/2) \in\mathbb{Q}$$ This is untrue, so there are no periodic points in this case. \\

\begin{figure}
  
    \centering
    \hbox{\hspace{-.6in}\includegraphics[scale = .3]{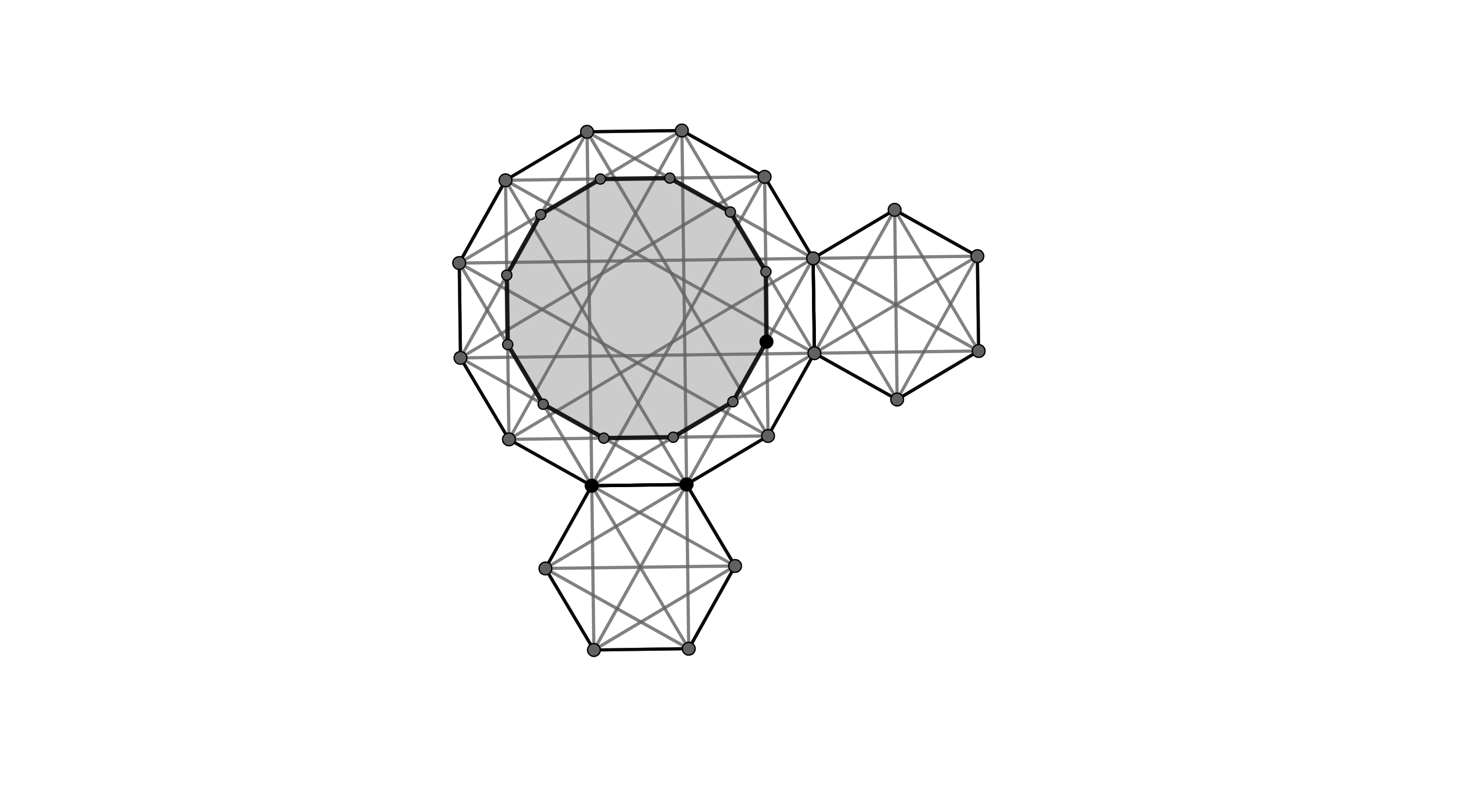}}
    \caption{The points in the shaded region of the surface of $(Y_{3,6},\eta_{3,6})$ are ruled out as possible periodic points, with the exception of the center}
    \label{fig:my_label9}
\end{figure}

\noindent Indeed, we examine Figure 9 and see case 1a applies to the most interior dodecagon of the figure, at least one the adjacent triangles, and at least one of the kites adjacent to those triangles. Case 1b corresponds to at least one of the kites adjacent to the kite of case 1a. Case 2a corresponds to at least one of the irregular hexagons inside shaded region, and case 2b corresponds to a triangle adjacent to those hexagons in the shaded region. By applying $\psi$, therefore, all points in the shaded region except the origin are ruled out as periodic. \\

\noindent Notice that this region covers exactly half of the horizontal length of $C_1$ near the $x$-axis and more than $1/3$ of the horizontal distance at vertical distances of over $1/3$ away from the axis. It follows from lemmas 2.4, 2.5, and 2.6 that the $\phi$-orbit of any point of $C_1$ intersects this doubly-shaded region. Therefore, the only points of $C_1$ away from the boundaries that are potentially periodic have $p_x = 0$. But, as already stated, half of the length of the locus $p_x = 0$ (except the origin) has been ruled out as periodic, and its image under $\phi$ is the other half (minus the center of the hexagon $C_1$ intersects in its interior). So, the only (potentially) periodic points in this locus are the origin and the center of the hexagon. The images of $C_1$ under $\langle\psi\rangle$ cover $Y_{3,6}$ away from the singularities, meaning no points but the singularities and the centers of the polygons are periodic under the action of the affine group.

\end{proof} 

\section{The Case of Odd $n$} 

\noindent Without much modification, the proof of theorem 5.1 can be generalized to prove theorem 1.4. Again, the proof of the theorem will be based on the pattern of cylinders that $\psi^k(p)$ traverses as $k$ varies - or rather - the pattern of restrictions on the coordinates of $\psi^k(p)$ obtained by applying the rational height lemma. Roughly, the proof will be divided into three steps: first, reducing to the cases of a few simple height-rule patterns, and second, ruling out points of the $2n$-gon following those patterns as periodic, and third, making certain every point on the surface has an image into the locus studied in the second part of the proof.\\

\begin{proof}[Proof of Theorem 1.4] Translate and scale $(\ye)$ so that the origin of the plane is at the center of the $2n$-gon, and so that the $2n$-gon is inscribed in the unit circle. Note that the interior angle measures of the $2n$-gon are $(2n-2)\pi/(2n)$. Now, let $T$ be the triangle with vertices at the origin, $(1,0)$, and $(\cos(-\pi/n),\sin(-\pi/n))$. As noted previously, this surface admits a horizontal cylinder decomposition. Notice that because $n$ is odd, two of the cylinders have their boundaries on the $x$-axis. Call the union of these cylinders $L_1$. Note $T\subset L_1$. Furthermore, these two cylinders both have height $\sin(\pi/n)$. Therefore, any point $p\in L_1$ obeys $p_y/\sin(\pi/n)\in\mathbb{Q}$. Also, note that any point $p$ of the $2n$-gon has an image $p^\prime = \psi^k(p) \in T$. It follows that $p^\prime,\psi(p^\prime)\in L_1$. \\

\noindent Now, we will make observations about the possible patterns of cylinders traversed by $p\in T$ under the orbit of $\langle \psi\rangle$. Let $L_2$ be the cylinder lying above $L_1$ in the $2n$-gon, and let $L_3$ be the cylinder lying below $L_1$ in the $2n$-gon. See Figure 10 for an example. Let $p \in T$. Note that either $\psi^{-1}(p)\in L_1$ or $L_3$, $\psi(p)\in L_1$, and $\psi^2(p)\in L_1$ or $L_2$. \\

\begin{figure}[t]
  
    \centering
    \hbox{\hspace{-.5in}\begin{overpic}[scale = .6]{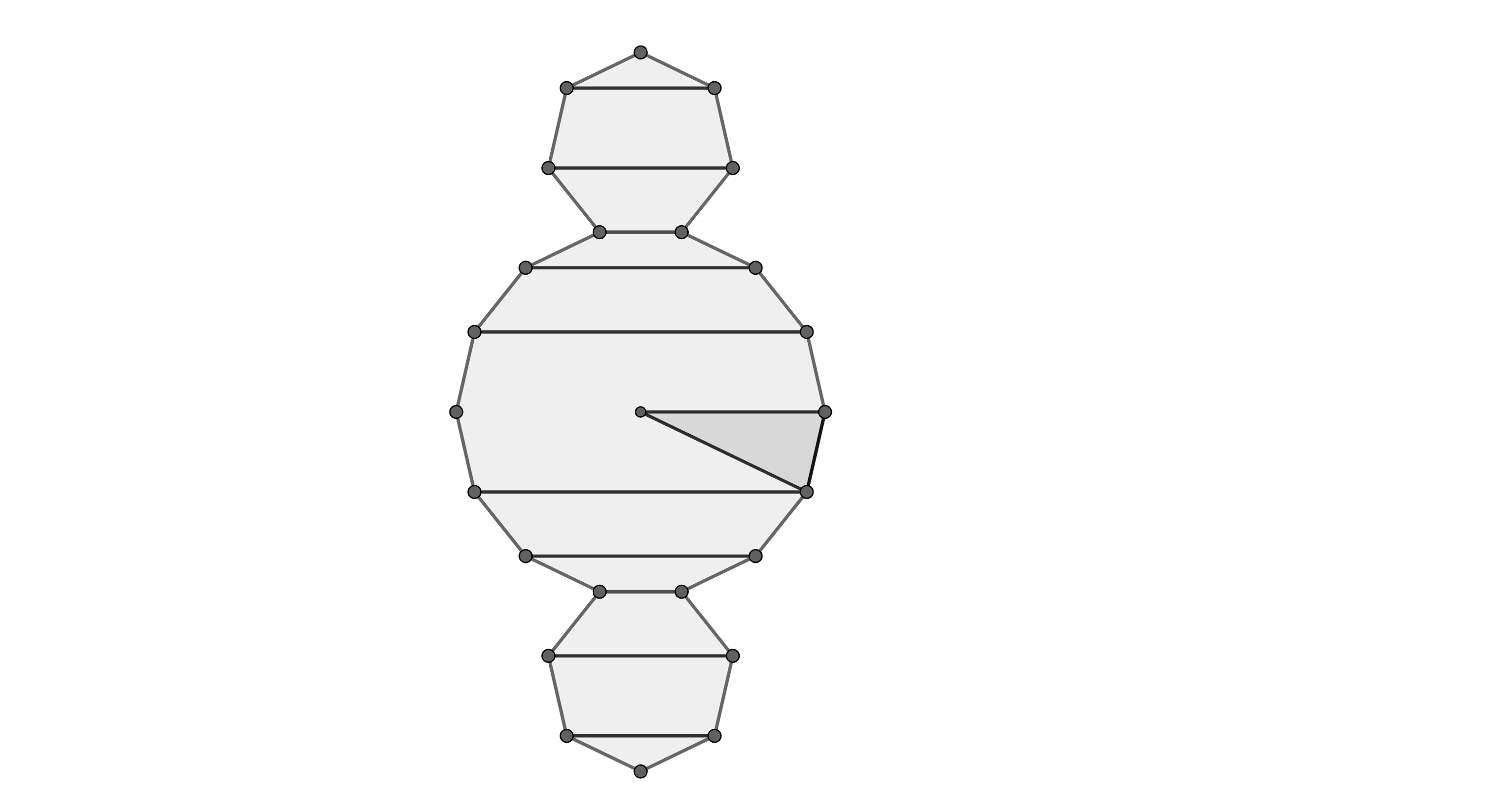}
    \put (36,34) {$L_2$}
    \put (46,29) {$L_1$}
    \put (52,24) {$T$}
    \put (35,19) {$L_3$}
    \end{overpic}}
    \caption{The regions $L_1, L_2, L_3$ and $T$ in the surface $(Y_{3,7},\eta_{3,7})$}
\end{figure}
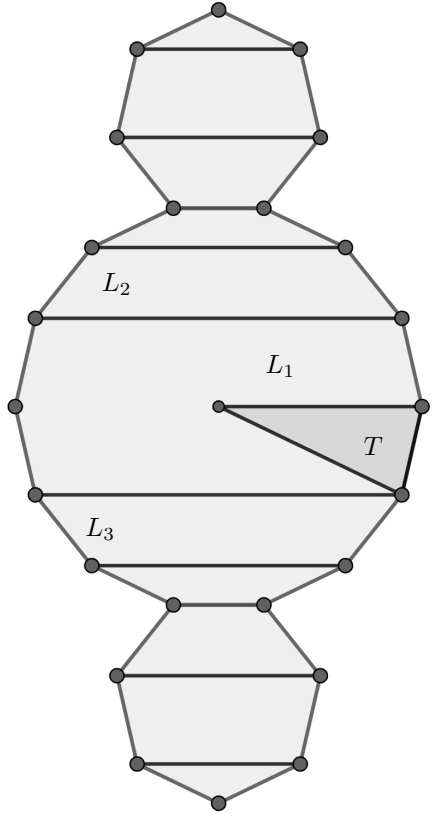

\noindent Now, let $p\in T$ and suppose $p$ is periodic. It follows $p_y/\sin(\pi/n) = q \in\mathbb{Q}$, and because $\psi(p) \in L_1$, $\psi(p)_y = (\sin(\pi/n)p_x + \cos(\pi/n)p_y )/\sin(\pi/n) = r \in \mathbb{Q}$. In other words, $p_x = r - q\cos(\pi/n) $. Based on our observations about the $\langle\psi\rangle$ orbit of $p$, we break our analysis into cases.\\

\noindent Case 1: Suppose $\psi^{-1}(p)\in L_1$. Then, $(-\sin(\pi/n)p_x + \cos(\pi/n)p_y )/\sin(\pi/n) = s \in \mathbb{Q}$ by the rational height lemma. In other words, $p_x = s - q \cos(\pi/n)$. Subtracting with the equation $p_x = r - q\cos(\pi/n) $ gives $$0 = s-r -2q \cos(\pi/n)$$ Because $n\geq 5$, $\cos(\pi/n)$ is irrational, and therefore $q = p_y = 0$. It follows from the definition that $p_x = r$.\\

\noindent Case 1a: Suppose that $\psi^2(p)\in L_1$. Then, $\psi^2(p)_y = 2\sin(\pi/n)\cos(\pi/n)p_x$ must satisfy $\psi^2(p)_y/\sin(\pi/n) \in \mathbb{Q}$, so $2r\cos(\pi/n) \in \mathbb{Q}$. Again, because cosine is irrational, this means $r = 0$, so the only periodic point in this case is $p_x = p_y = 0$, the origin.\\

\noindent Case 1b: Suppose that $\psi^2(p)\in L_2$. Then, by the rational height lemma, $$ \frac{2\sin(\pi/n)\cos(\pi/n)p_x -\sin(\pi/n)}{2\sin(\pi/n)\cos(\pi/n)-\sin(\pi/n)} = \frac{2\cos(\pi/n)p_x - 1}{2\cos(\pi/n) -1} \in \mathbb{Q}$$ This is only possible if $p_x = 1$, a singularity. \\

\noindent Case 2: Suppose that $\psi^{-1}(p)\in L_3$ and that $\psi^2(p)\in L_1$. The height rule applied to this second fact gives $$\frac{2\sin(\pi/n)\cos(\pi/n) p_x + (\cos^2(\pi/n)-\sin^2(\pi/n))p_y}{\sin(\pi/n)} $$ $$ = 2\cos(\pi/n)r - 2\cos^2(\pi/n) q -(\cos^2(\pi/n)-\sin^2(\pi/n))q  $$ $$ = 2\cos(\pi/n)r - q \in\mathbb{Q}$$ It follows that $r = 0$. Now, from $\psi^{-1}(p)\in L_3$, one has from the rational height lemma that $$\frac{-\sin(\pi/n)p_x + \cos(\pi/n)p_y + \sin(\pi/n)}{2\sin(\pi/n)\cos(\pi/n) -\sin(\pi/n)}= \frac{2\cos(\pi/n) q + 1}{2\cos(\pi/n) -1 } \in\mathbb{Q}$$ It follows that $q = -1$. Now, $r = 0$, $q = -1$ corresponds to the point $(p_x,p_y) = (\cos(\pi/n),-\sin(\pi/n))$, another singularity. \\

\noindent Case 3: Suppose that $\psi^{-1}(p)\in L_3$ and that $\psi^2(p)\in L_2$. According to the rational height lemma, we can obtain $u,v\in\mathbb{Q}$ as $$ u = \frac{-\sin(\pi/n)p_x+ \cos(\pi/n)p_y + \sin(\pi/n)}{2\sin(\pi/n)\cos(\pi/n) -\sin(\pi/n)} = \frac{2\cos(\pi/n)q -r + 1}{2\cos(\pi/n) -1} $$ and $$v =\frac{2\sin(\pi/n)\cos(\pi/n)(r-\cos(\pi/n)q) + (\cos^2(\pi/n) - \sin^2(\pi/n))q -\sin(\pi/n)}{2\sin(\pi/n)\cos(\pi/n) -\sin(\pi/n)}$$ $$ = \frac{2\cos(\pi/n)r -q -1}{2\cos(\pi/n) -1} $$ Examine $u$, for example. It must be that $u = q$, in which case $-u = -q = -r + 1 $, so $r-q = 1$. This gives a linear relation between $p_x$ and $p_y$, namely: $$\frac{\sin(\pi/n)}{1-\cos(\pi/n)}(p_x - 1) = p_y$$ In other words, $p$ must lie on the boundary of the $2n$-gon.\\

\noindent In summary, all interior points of the $2n$-gon except the origin have been ruled out as potential periodic points. Lastly, observe that by lemmas 2.4, 2.5, and 2.6, every point of $Y_{3,n}$ except for the singularities has image into the interior of the $2n$-gon away from the origin. Therefore, only the origin and singularities are periodic.

\end{proof}

\noindent For corollaries we have:

\begin{corollary} There exists an infinite sequence of primitive Veech surfaces $(X_i,\omega_i)$ such that the genus of $X_i$ is unbounded and the number of periodic points of the $(X_i,\omega_i)$ is uniformly bounded.

\end{corollary}

\begin{proof} Letting $(X_i,\omega_i) = (Y_{3,2i+1},\eta_{3,2i+1})$, theorem 1.4 guarantees each $(X_i,\omega_i)$ has at most $1$ periodic point which is not a singularity. Note that if $3|(2i+1)$, $(X_i,\omega_i)$ has three singularities, and otherwise, it has only one. In summary, the $(X_i,\omega_i)$ have at most $4$ periodic points each. Also, the genus of $(X_i,\omega_i)$ is $$\frac{ 3(2i+1) - (2i + 1) - 3 - (2i+1,3) }{2}  + 1 \geq 2i - 1$$ which is unbounded. \\

\noindent Note the surfaces $(\ye)$ are primitive by theorem 3.3, completing the proof of the corollary.

\end{proof}

\noindent Also: \begin{corollary} There exists an infinite sequence of primitive Veech surfaces $(X_i,\omega_i)$ such that the genus of $X_i$ is unbounded in $i$ and such that not all Weierstra{\ss} points are periodic.

\end{corollary}

\begin{proof} The genus of $ (Y_{3,2i+1},\eta_{3,2i+1})$ is $$\frac{ 3(2i+1) - (2i + 1) - 3 - (2i+1,3) }{2}  + 1 \geq 2i - 1$$ and therefore has at least $4i = 2(2i-1) + 2$ Weierstra{\ss} points. On the other hand, $ (Y_{3,2i+1},\eta_{3,2i+1})$ has at most $4$ periodic points by theorem 1.4, noting that such a surface has at most $3$ singularities.

\end{proof}

\noindent The versions of these corollaries without the ``primitive" qualifier are not very interesting. Indeed, there is (at worst) a degree-$2$ branched cover of any translation surface that admits cyclic branched covers, branched exactly over the periodic points, mapping this locus bijectively to the coverings. In particular, these translation coverings (at worst) double the number of periodic points from the underlying surface to the covers but are unbounded in complexity.

\section{The Case of Even $n$, $n\geq 8$} 

\noindent Let $n\geq 8$. Translate and scale $(\ye)$ so that its sides have length $1$ and its center is a the origin. Again, note the interior angle measures of the $2n$-gon are $(2n-2)\pi/(2n)$. Now, because $n$ is even, the $2n$-gon of $(\ye)$ has two vertical sides, and one of the $n$-gons may be viewed as having two vertical sides being glued on to these. These vertical sides form a height of a cylinder $C_0$, which contains the origin. The widths of the $2n$-gon and the $n$-gons are the same as their heights, which are given by $\cot(\pi/2n)$ and $\cot(\pi/n)$ respectively. (To see this, take the right triangle with one vertex at the center of a polygon, one vertex at the midpoint of a vertical side of the polygon, and one vertex at an endpoint of the same side of the polygon, and compute the tangent of the angle at the center of the polygon.) Consequently, the length of the cylinder $C_0$ is $$\cot(\frac{\pi}{2n}) + \cot(\frac{\pi}{n}) \leq \frac{1}{\sin(\frac{\pi}{2n})}  + \frac{1}{\sin(\frac{\pi}{n})} <  \frac{3}{\sin(\frac{\pi}{n})}$$ with the last inequality deriving from the fact $\sin(x/2) > \sin(x)/2$ for positive $x$ near $0$. Let $T\subset C_0$ be the triangle with vertices $(\cot(\frac{\pi}{2n})/2,-1/2),(0,0)$ and $(\cot(\frac{\pi}{2n})/2,1/2)$. In other words, the union of the images of $T$ under $\psi$ is the entire $2n$-gon of $(\ye)$. \\

\noindent Let the cylinders of the decomposition be numbered so that traveling vertically through the $2n$-gon from the base, crossing a cylinder boundary increases the index of the cylinder $C_i$ by $1$. In other words, in this numbering, the cylinders above and below $C_0$ in the $2n$-gon are $C_1$ and $C_{-1}$ respectively, and they have heights $\sin(\pi/2 - \pi/n) = \cos(\pi/n)$. Above $C_1$ in the $2n$-gon is $C_2$ with height $\sin(\pi/2 - 2\pi/n) = \cos(2\pi/n)$, and below $C_{-1}$ is $C_{-2}$ also with height $\cos(2\pi/n)$. Note that the base of $C_2$ is at height $1/2 +\cos(\pi/n) $, and the boundary of $C_{-2}$ begins at $-1/2 -\cos(\pi/n)$. Further, note that the images of a point $p\in T$ under consecutive powers of $\psi$ must be in the same cylinder or cylinders adjacent in the $2n$-gon. \\

\noindent We are now ready to show: \begin{theorem}
Let $n\geq 8$ be even. The Ward surface $(\ye)$ has either $6$ or $4$ periodic points, depending on whether $3|n$ or not, respectively. The periodic points are the vertices together with the centers of the polygons in the regular polygonal decomposition.
\end{theorem}

\begin{proof} Let $p\in T$. The proof will start by ruling out $p$ as periodic depending on the pattern of cylinders $\psi^k(p)$ lie in for varying $k$. After a sufficient number of cases have been examined, we will have a locus $L$ of non-periodic points whose images under the affine group cover the surface - minus vertices and the centers of the polygons, completing the theorem. Assume (for contradiction) now that $p$ is periodic. Let us break into a case-by-case study, and in any case, $p\in T\subset C_0$, so $p_y = q\in\mathbb{Q}$ by the rational height lemma. Notice that $\psi(p)\in C_0$ or $C_1$ as $p\in T$. These are our two main cases.\\

\noindent Case 1: Suppose $\psi(p)\in C_0$. Note for later then that $\psi^2(p)\in C_0$ or $C_1$ and, since $p\in T$, that $\psi^{-1}(p) \in C_{-1}$ or $C_0$. This gives us four sub-cases to analyze. In any sub-case, as $\psi(p)\in C_0$, the rational height lemma says $\psi(p)_y = \sin(\pi/n)p_x + \cos(\pi/n)p_y = r\in\mathbb{Q}$. So, $\sin(\pi/n)p_x = r- \cos(\pi/n)q$.\\

\noindent Case 1a: Suppose further that $\psi^{-1}(p), \psi^2(p)\in C_0$. Then, by the rational height lemma, $$-\sin(\pi/n)p_x + \cos(\pi/n)p_y \in \mathbb{Q}$$ and $\sin(2\pi/n)p_x + \cos(2\pi/n)p_y\in\mathbb{Q}$. Adding $$\sin(\pi/n)p_x + \cos(\pi/n)p_y$$ and $$-\sin(\pi/n)p_x + \cos(\pi/n)p_y$$ it follows that $2\cos(\pi/n)p_y\in\mathbb{Q}$, but $p_y = q\in\mathbb{Q}$, hence $p_y = 0$ by the irrationality of the cosine. We are left with the facts that $$\sin(\pi/n)p_x,2\sin(\pi/n)\cos(\pi/n)p_x\in\mathbb{Q}$$ If $p_x\neq 0$, then dividing yields $2\cos(\pi/n)\in\mathbb{Q}$, which is false. Therefore, $p_x,p_y = 0$, so the origin is the only possible periodic point in this case.\\

\noindent Case 1b: Suppose that $\psi^{-1}(p)\in C_0$, $\psi^2(p)\in C_1$. Again, $$-\sin(\pi/n)p_x + \cos(\pi/n)p_y \in \mathbb{Q}$$ by the rational height lemma applied to $\psi^{-1}(p)$, and adding $\sin(\pi/n)p_x + \cos(\pi/n)p_y$ and $-\sin(\pi/n)p_x + \cos(\pi/n)p_y$, it follows that $2\cos(\pi/n)p_y\in\mathbb{Q}$, but $p_y = q\in\mathbb{Q}$, hence $p_y = 0$ by the irrationality of the cosine. Now, knowing that $p_y =q = 0$, the rational height lemma applied to $\psi^2(p)$ yields $$(2\sin(\pi/n)\cos(\pi/n)p_x - 1/2)/\cos(\pi/n)  = (2r\cos(\pi/n) - 1/2)/\cos(\pi/n)\in\mathbb{Q}$$ This is impossible, and therefore there are no periodic points in this case.\\

\noindent Case 1c: Suppose that $\psi^2(p)\in C_0$ but $\psi^{-1}(p)\in C_{-1}$. Apply either Case 1a or Case 1b to $\psi(p)$.\\

\noindent Case 1d: Suppose $\psi^{-1}(p)\in C_{-1} $ and $\psi^2(p)\in C_1$. The rational height lemma applied to $\psi^2(p)$ yields $$\mathbb{Q}\ni (2\sin(\pi/n)\cos(\pi/n)p_x  + (\cos(\pi/n)^2 - \sin(\pi/n)^2)p_y- 1/2)/\cos(\pi/n)$$ which by the Pythagorean theorem is the same as $$\frac{2\cos(\pi/n) r - q - 1/2}{\cos(\pi/n)}$$ By the irrationality of the cosine, this is impossible unless $q = -1/2$, which only happens at the vertex in $T$. \\

\noindent Case 2: Suppose $\psi(p)\in C_1$. Then, by the rational height lemma $$\frac{\sin(\pi/n)p_x + \cos(\pi/n)p_y-1/2}{\cos(\pi/n)} = r\in\mathbb{Q}$$ or $$\sin(\pi/n)p_x = 1/2-\cos(\pi/n)q + \cos(\pi/n)r$$ First, in case 2a, we will rule out periodic points (except for a vertex) in the case $\psi^2(p)\in C_1$. It will follow that $\psi^2(p)\in C_2$. Note that because $p\in T$, $\psi^{-1}(p)\in C_0$ or $C_{-1}$. However, we omit the cases where $\psi^{-1}(p)\in C_0$ as we can apply $\psi$ to $p$ and land back in case 1. Therefore, the two cases left to study are $\psi^{-2}(p)\in C_{-1}$ or $\psi^{-2}\in C_{-2}$ with $\psi^{-1}(p)\in C_{-1}$ in either case. These last two cases are taken up respectively in case 2b and the subsequent paragraphs.\\

\noindent Case 2a: Suppose further that $\psi^2(p)\in C_1$. Then, $$(2\sin(\pi/n)\cos(\pi/n)p_x  + (\cos(\pi/n)^2 - \sin(\pi/n)^2)p_y- 1/2)/\cos(\pi/n) = t\in\mathbb{Q}$$ Using that $$\sin(\pi/n)p_x = 1/2-\cos(\pi/n)q + \cos(\pi/n)r$$ and utilizing the Pythagorean theorem, then multiplying by cosine, $$t\cos(\pi/n) = -q -1/2+ \cos(\pi/n) + 2\cos(\pi/n)^2r $$ or $$-q -1/2+ \cos(\pi/n)(1-t) + \cos(\pi/n)^2r = 0$$ However, because $n\geq 8$ the cosine does not satisfy a degree $2$ polynomial with rational coefficients, we must have that $r = 0$, $q = -1/2$ and $t = 1$. The reader is welcomed to verify this is a vertex, the only possible periodic point in this case.\\

\noindent Case 2b: Suppose $\psi^{-2}(p),\psi^{-1}(p)\in C_{-1}$. Then, by the rational height lemma applied to $\psi^{-1}(p)$ $$\frac{-\sin(\pi/n)p_x + \cos(\pi/n)p_y-1/2}{\cos(\pi/n)} = t\in\mathbb{Q}$$ or $$-\sin(\pi/n)p_x = 1/2-\cos(\pi/n)q + \cos(\pi/n)t$$ The rational height lemma applied to $\psi^{-2}(p)$ gives that $$(-2\sin(\pi/n)\cos(\pi/n)p_x  + (\cos(\pi/n)^2 - \sin(\pi/n)^2)p_y- 1/2)/\cos(\pi/n) = u\in\mathbb{Q}$$ Utilizing the Pythagorean theorem, then multiplying by cosine, $$u\cos(\pi/n) = -q -1/2+ \cos(\pi/n) + 2\cos(\pi/n)^2t$$ Like above, this is not possible  unless $q = -1/2$, $u = 1$, and $t = 0$, again a vertex.\\

\noindent In summary, we have ruled out every point of $T$ being periodic except the origin and points $p\in T$ having the property that $\psi^2(p)\in C_2$ and $\psi^{-2}(p) \in C_{-2}$. Let us study such a point $p$. Well, then, it  must be the case that $$\sin(2\pi/n)p_x + \cos(2\pi/n)p_y \geq \cos(\pi/n) + 1/2 $$ and $$-\sin(2\pi/n)p_x + \cos(2\pi/n)p_y \leq -\cos(\pi/n) - 1/2 $$ Multiplying the second expression by $-1$, adding and dividing by $2$ yields $$\sin(2\pi/n)p_x \geq \cos(\pi/n) + 1/2 $$ Using the double angle formula, we see $$\sin(\pi/n)p_x \geq 1/2 + 1/(4\cos(\pi/n)) \geq 3/4$$ Or $$p_x\geq \frac{3}{4\sin(\pi/n)}$$ It follows that all points of $T$ with $0 \leq  p_x \leq  \frac{3}{4\sin(\pi/n)}$  except the origin have been ruled out as being periodic. Call this sub-locus $S$, and let $L$ be the union of all images of $S$ under $\psi$. See Figure 11 for an example.\\

\begin{figure}[t]
  
    \centering
    \hbox{\hspace{.5in}\begin{overpic}[scale = .23]{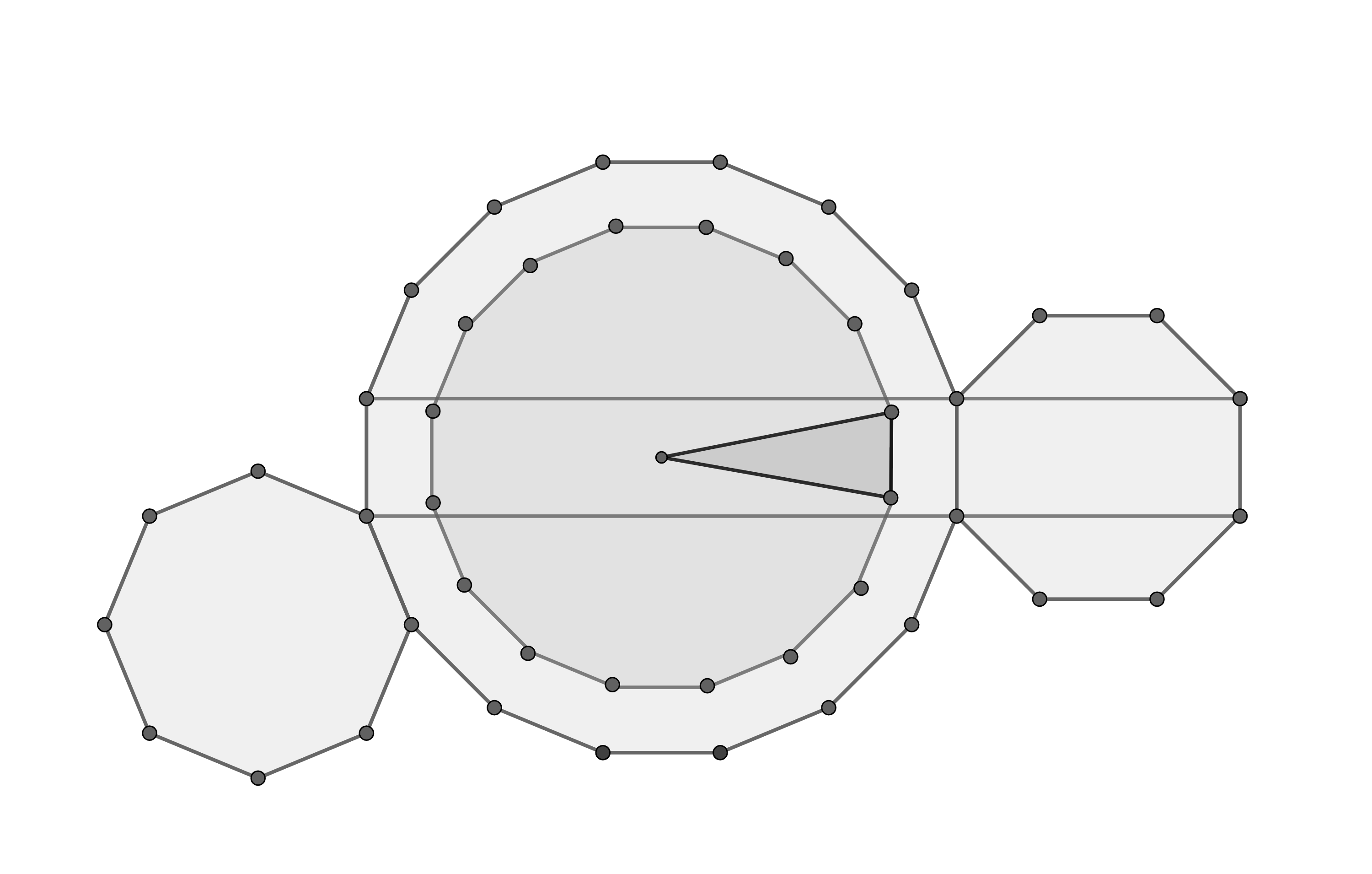}
    \put (29,36) {$M$}
    \put (80,29) {$C_0$}
    \put (62,29) {$S$}
    \put (37,22) {$L$}
    \put (65,36) {$B$}

    \end{overpic}}
    \caption{The surface $(Y_{3,8},\eta_{3,8})$ labeled with regions introduced in the proof of theorem 8.1}
\end{figure}
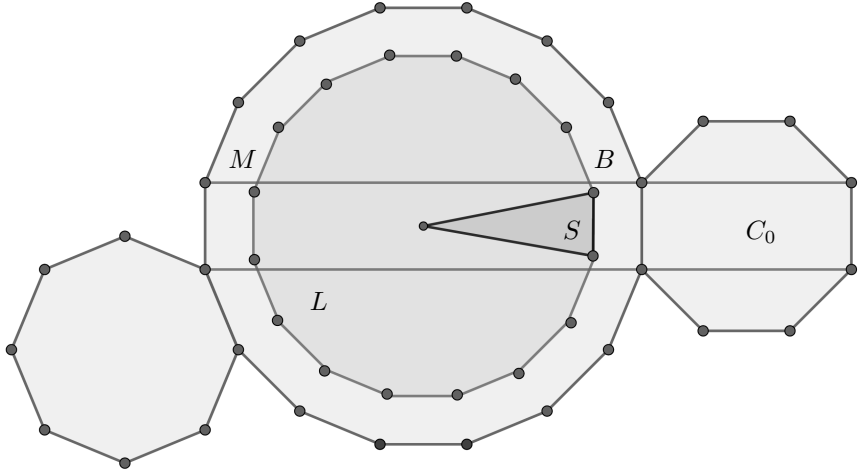

\noindent This polygon $L$ is similar to the $2n$-gon of $(\ye)$ centered at the origin. Some simple estimates will show that every point $p^\prime$ of $C_0$ away from a boundary has an image by a power of $\phi$, the reducible affine map with horizontal fixed direction, into $L$, away from the origin as long as $p^\prime$ is not the center of the $n$-gon, and therefore that the only potentially periodic points of $C_0$ away from the boundaries are the centers of the polygons. The estimates will show that at vertical heights between $-1/3$ and $1/3$ in $C_0$, a leaf of the horizontal foliation has more than half of its length in $L$, and at vertical heights between $1/3$ and $1/2$ and $-1/2$ and $-1/3$, a leaf of the horizontal foliation has more than a third of its length in $L$. By lemmas 2.4, 2.5, and 2.6, these estimates will show that the only possible periodic points of $C_0$ are the centers of the polygons and the vertices. Moreover, because the set of the centers of the polygons is preserved by the affine group, these points are in fact periodic. Also, the union of the images of the interior of $C_0$ under powers of $\psi$ cover the surface $(\ye)$, meaning the only periodic points are the vertices and the centers of the polygons. \\

\noindent So, it remains to demonstrate that any point of $C_0$ away from the boundary has an image under a power of $\phi$ into $L$, away from the origin if the point is not a center of the $n$-gon. $S$ is the subset of $T$ such that $0\leq p_x \leq \frac{3}{4\sin(\pi/n)}$, which is similar to $T$. It follows that $L$ has width $$ 2( \frac{3}{4\sin(\pi/n)}) = \frac{3}{2\sin(\pi/n)},$$ which is strictly larger than half the length of $C_0$, as our computations made just before the proof demonstrate.\\

\noindent We want a lower bound on the length of the vertical side of $S$ that does not depend on $n$. We can find one by studying the scaling factor of similarity between $L$ and the $2n$-gon. The $2n$-gon has horizontal length $$\cot(\frac{\pi}{2n}) < 2/\sin(\pi/n)$$ and $L$ has width  $\frac{3}{2\sin(\pi/n)}$. Therefore, the ratio between the horizontal lengths of $L$ and the $2n$-gon is at least $$ \frac{3}{2\sin(\pi/n)} /\cot(\frac{\pi}{2n}) > \frac{3}{2\sin(\pi/n)} \frac{\sin(\pi/n)}{2} = 3/4.  $$ Since the vertical height of $T$ is $1$, the vertical height of $S$ is at least $\frac{3}{4}$, and so every point with $|p_y| \leq \frac{3}{8}$ has image by a power of $\phi$ into $L$. This image is away from the origin if $p$ is not the origin itself or the center of the $n$-gon.\\

\noindent Next, we consider the locus where $ 3/8 \leq p_y < 1/2$. Because $1/3 < 3/8$, we simply have to show that the horizontal length of $L$ at a given height in this locus is at least $1/3$ the length of $C_0$, and then lemmas 2.4, 2.5, 2.6 will show any such point in this locus $ 3/8 \leq p_y < 1/2$  in $(\ye)$ has image by a power of $\phi$ into $L$. The problem is that the horizontal length of $L$ is potentially decreasing as $p_y$ increases from $ 3/8 \leq p_y < 1/2$, but luckily, it does not decrease by too much. Let $B$ denote the side of $L$ which is the image of the vertical side of $S$ by $\psi$. $B$ is similar to the image of the vertical side of $T$ under $\psi$, so therefore, $B$ has change in height at least $$\frac{3}{4}\cos(\pi/n) \geq \frac{3}{4}\cos(\pi/8) > 5/8$$ Therefore, the higher endpoint of $B$ lies above the boundary of $C_0$. Let $M$ be $B$ reflected across the $y$-axis. We are interested in the ratio of the horizontal distance between the upper endpoints of $M$ and $B$ to the horizontal distance between the lower endpoints of $M$ and $B$. Because $n\geq 8$, $M$ and $B$ are joined above by a path of at least $(2n-2)/2 -2 = n-3 \geq 5$ other sides of $L$, all having horizontal length longer than $M$ and $B$. If $\Sigma$ is the sum of horizontal lengths of segments in the path, and $l(B)$ is the horizontal length of $B$, then the ratio of lengths is $$\frac{\Sigma}{2l(B) + \Sigma} = \frac{1}{\frac{2l(B)}{\Sigma} + 1}$$ and since $$\frac{2l(B)}{\Sigma} \leq \frac{2l(B)}{(n-3)l(B)}$$ we have: $$\frac{1}{\frac{2l(B)}{\Sigma} + 1} \geq \frac{1}{\frac{2l(B)}{(n-3)l(B)} + 1} \geq \frac{5}{7}$$ Therefore, the ratio we are looking for is at least $5/7$. \\

\noindent In summary, the horizontal length of $L$ at a given height between $3/8$ and $1/2$ is at least $(\frac{1}{2})(\frac{5}{7}) $ times the length of $C_0$ which is greater than $1/3$ by $1/42$.\\

\noindent The locus $-1/2 < p_y \leq -3/8$ may be dealt with similarly, or by applying $\psi^n$ and the above paragraphs again.

\end{proof}

\section{Applications}
\subsection{Applications to Sections of Bundles Constructed from Veech Groups}
\noindent In \cite{Moeller04}, M\"{o}ller has shown that periodic points of Veech surfaces correspond with holomorphic sections of certain surface bundles over the appropriate Teichm\"{u}ller curve. In this way, our computations bound the number of sections of certain bundles whose constructions we now describe. The exposition will follow that of Shinomiya in \cite{syarticle}.\\

\noindent Let $(X,\omega)$ be any translation surface. Mark $X$ by a homeomorphism $\iota:S_g\mapsto X$. This defines a map $F:\mathbb{H}^2 \mapsto T(S_g)$, the Teichm\"{u}ller space of a genus $g$ surface as follows. View the upper half plane $\mathbb{H}^2 = SL_2(\mathbb{R})/SO(2)$. For $a\in \mathbb{H}^2$, choose a representative $A\in SL_2(\mathbb{R})$. Set $F(a) = (i,AX)\in T(S_g)$ where $AX$ is the Riemann surface underlying the translation surface $A(X,\omega)$. This definition is independent of choices, and it is well known that $F$ is a holomorphic isometry onto its image called a \textit{Teichm\"{u}ller disk}. Note that taking the quotient by the mapping class group yields a map $$f_0: \mathbb{H}^2/PSL(X,\omega) \mapsto M_g$$ where $PSL(X,\omega)$ is the quotient of $SL(X,\omega)$ by $\langle -I\rangle$ and $M_g$ is the moduli space of curves of genus $g$. Notice that $PSL(X,\omega)$ is discrete by a result of Veech.\\

\noindent Now, let $(X,\omega)$ be a Veech surface. Because $(X,\omega)$ is a Veech surface, $PSL(X,\omega)$ is finitely generated and therefore has a finite index subgroup $\Gamma$ without torsion. Set $B = \mathbb{H}^2/\Gamma$.\\

\noindent There is a map $f: B\mapsto M_g$ defined by taking the quotient by $PSL(X,\omega)$ and composing with $f_0$. In other words, associated to each $b\in B$ is a Riemann surface $AX$, where $A$ depends on $b$. One translates this data into $$M = \{(b,p)\,:\, b\in B, p\in f(b) = AX\}$$ a holomorphic family of Riemann surfaces with base space $B$ (and call the projection $\pi$). A holomorphic family of Riemann surfaces constructed in this way is said to be \textit{constructed from a Veech group}. Notice that for $a_0 = [i] \in B$, $\pi^{-1}(a_0) = X$ (where $i= \sqrt{-1}$ imagined in $\mathbb{H}^2$). In other words, a section evaluated at $a_0$ yields a point in our original surface. \\

\noindent Let $S$ stand for the set of all holomorphic sections of the bundle $(M,\pi,B)$, and let $s\in S$. Shinomiya (see \cite{syarticle}) has the following theorem:

\begin{theorem}[Shinomiya] For a bundle $(M,\pi,B)$ constructed from a Veech surface $(X,\omega)$ (as above), the map $s\mapsto s(a_0)$ is injective. Moreover, the point $s(a_0)$ is a periodic point of $(X,\omega)$.

\end{theorem}

\noindent This allows us to restate many of our theorems and corollaries as statements about the number of sections of certain bundles. For example:

\begin{corollary} Let $n \geq 5$ be odd. From the Riemann surface $Y_{3,n}$ defined by the equation $$y^{2n} = (u-2)(u+1)^2$$ with translation structure given by the differential $$(u+1)du/y^{2n-1}$$ construct a holomorphic family of surfaces from the Veech group. This bundle has at most four or two sections, depending on if $3|n$ or not, respectively.

\end{corollary}

\begin{proof} From theorems 1.4 and 9.1, we know this bundle has at most 4 or 2 sections depending on whether $3|n$ or not, respectively.

\end{proof}

\subsection{Applications to the Study of Infinitely Generated Veech Groups} 

\noindent In \cite{Hubert2004INFINITELYGV}, Hubert and Schmidt give a construction that yields translation surfaces with infinitely generated Veech groups. To introduce the construction, we need the notion of a \textit{connection point} on a translation surface: A point on a translation surface is a \textit{connection point} if every separatrix passing through it is a saddle connection.\\

\noindent We also need the notion of \textit{weak} and \textit{strong holonomy types} as defined in \cite{Hubert2004INFINITELYGV}: A translation surface is of weak holomony type if all holonomy vectors and saddle connection vectors, when expressed in terms of the canonical basis of $\mathbb{R}^2$, have all components in the holomony field of the translation surface. For the notions of holonomy field, holonomy vectors, and saddle connection vectors, see \cite{Hubert2004INFINITELYGV}. A translation surface is of \textit{strong holonomy type} if it is of weak holonomy type and the periodic directions are exactly the vertical direction and those directions whose slopes belong to the holonomy field.\\

\noindent We can now state the following theorem of Hubert and Schmidt.\\

\begin{theorem}[Hubert and Schmidt] Let $(X,\omega)$ be a nonarithmetic Veech surface of strong holonomy type. Then there are infinitely many nonperiodic connection points on $(X,\omega)$. Futher, if $p$ is one of these points, any translation covering of $(X,\omega)$ branched over $p$ (along with periodic points of $(X,\omega)$)  has an infinitely generated Veech group.

\end{theorem}

\noindent Hubert and Schmidt (in \cite{Hubert2004INFINITELYGV}) note a theorem of McMullen that gives a criterion for a translation surface to be of strong holonomy type:

\begin{theorem}[McMullen] A Veech surface having real quadratic holonomy field lies in the $GL_2(\mathbb{R})$-orbit of a nonarithmetic Veech surface of strong holonomy type. \end{theorem}

\noindent We have the following corollary:

\begin{corollary}If $n = 4, 5$ or $6$, the Ward surfaces $(\ye)$ are in the $GL_2(\mathbb{R})$-orbit of a surface of strong holonomy type.
\end{corollary} 

\begin{proof} Hooper in \cite{hooper} notes that the holonomy field of the surface $(\ye)$ is $\mathbb{Q}(\cos(\pi/n))$, which is real quadratic when $n = 4,5$ or $6$. Here, we are also using the fact that when the translation surface admits an affine diffeomorphism with hyperbolic derivative, the trace of the derivative generates the holonomy field - see \cite{Kenyon2000BilliardsOR}. Therefore, we may apply theorem 9.3.

\end{proof}

\noindent Another theorem of Hubert and Schmidt from \cite{Hubert2004INFINITELYGV} makes it simple to identify connection points on a Veech surface of strong holonomy type. Recall that a \textit{rational point} of a translation surface is a point that has finite orbit under each of a pair of affine maps with parabolic derivatives having transverse fixed directions. Now consider:

\begin{theorem}[Hubert and Schmidt] Let $p\in (X,\omega)$ be a nonsingular point on a Veech surface of strong holonomy type. Then $p$ is a connection point if and only if $p$ is a rational point. \end{theorem}

\noindent In particular, it is easy to find connection points on $(\yef),(\yefi),$ and $(Y_{3,6},\eta_{3,6})$ because it is easy to find rational points. Because we have classified the periodic points of these surfaces, it is easy to find non-periodic connection points and therefore construct new examples of translation surfaces with infinitely generated Veech groups.

\section{Acknowledgements}

\noindent The author would like to thank Chris Leininger for suggesting the problem of computing periodic points and Autumn Kent for careful readings of the early forms of this paper. Most of all, the author wishes to thank his friends and his parents for their patience and support.

\bibliographystyle{unsrtnat}


\begin{thebibliography}{14}
\providecommand{\natexlab}[1]{#1}
\providecommand{\url}[1]{\texttt{#1}}
\expandafter\ifx\csname urlstyle\endcsname\relax
  \providecommand{\doi}[1]{doi: #1}\else
  \providecommand{\doi}{doi: \begingroup \urlstyle{rm}\Url}\fi

\bibitem[Möller(2009)]{Moeller09}
Martin Möller.
\newblock Affine groups of flat surfaces.
\newblock \emph{Handbook of Teichmüller theory}, pages 369--387, 2009.
\newblock \doi{10.4171/055-1/11}.

\bibitem[Moeller(2004)]{Moeller04}
M.~Moeller.
\newblock Periodic points on {Veech} surfaces and the {Mordell}-{Weil} group
  over a {T}eichm{\"u}ller curve.
\newblock \emph{Inventiones mathematicae}, 165:\penalty0 633--649, 2004.

\bibitem[Gutkin et~al.(2003)Gutkin, Hubert, and
  Schmidt]{ASENS_2003_4_36_6_847_0}
Eugene Gutkin, Pascal Hubert, and Thomas~A. Schmidt.
\newblock Affine diffeomorphisms of translation surfaces : periodic points,
  fuchsian groups, and arithmeticity.
\newblock \emph{Annales scientifiques de l'\'Ecole Normale Sup\'erieure}, Ser.
  4, 36\penalty0 (6):\penalty0 847--866, 2003.
\newblock \doi{10.1016/j.ansens.2003.05.001}.
\newblock URL \url{www.numdam.org/item/ASENS_2003_4_36_6_847_0/}.

\bibitem[McMullen(2005)]{genus2veech}
Curtis McMullen.
\newblock Teichmüller curves in genus two: Discriminant and spin.
\newblock \emph{Mathematische Annalen}, 333, 09 2005.
\newblock \doi{10.1007/s00208-005-0666-y}.

\bibitem[Apisa et~al.(2020)Apisa, Saavedra, and Zhang]{apisa2020periodic}
Paul Apisa, Rafael~M. Saavedra, and Christopher Zhang.
\newblock Periodic points on the regular and double $n$-gon surfaces, 2020.

\bibitem[Phillip~Griffiths(1978)]{book:3970}
Joseph~Harris Phillip~Griffiths.
\newblock \emph{Principles of algebraic geometry}.
\newblock Pure and applied mathematics. Wiley, 1978.
\newblock ISBN 0471327921,9780471327929.
\newblock URL
  \url{http://gen.lib.rus.ec/book/index.php?md5=421ce3136b5da1dfaa45510e0f646047}.

\bibitem[Veech(1989)]{Veech89}
W.~A. Veech.
\newblock Teichm{\"u}ller curves in moduli space, eisenstein series and an
  application to triangular billiards.
\newblock \emph{Inventiones mathematicae}, 97:\penalty0 553--583, 1989.

\bibitem[Apisa(2020)]{aprhl20}
Paul Apisa.
\newblock Gl(2,r)-invariant measures in marked strata: Generic marked points,
  {Earle-Kra} for strata, and illumination.
\newblock \emph{Geometry \& Topology}, 24:\penalty0 373–408, 2020.
\newblock \doi{10.2140/gt.2020.24.373}.

\bibitem[Farkas and Kra(1980)]{nla.cat-vn1156072}
Hershel~M. Farkas and Irwin Kra.
\newblock \emph{Riemann surfaces}.
\newblock Springer-Verlag New York, 1980.
\newblock ISBN 0387904654.

\bibitem[Ward(1998)]{ward_1998}
Clayton~C. Ward.
\newblock Calculation of fuchsian groups associated to billiards in a rational
  triangle.
\newblock \emph{Ergodic Theory and Dynamical Systems}, 18\penalty0
  (4):\penalty0 1019–1042, 1998.
\newblock \doi{10.1017/S0143385798117479}.

\bibitem[Hooper(2012)]{hooper}
W.~Patrick Hooper.
\newblock {Grid Graphs and Lattice Surfaces}.
\newblock \emph{International Mathematics Research Notices}, 2013\penalty0
  (12):\penalty0 2657--2698, 05 2012.
\newblock ISSN 1073-7928.
\newblock \doi{10.1093/imrn/rns124}.
\newblock URL \url{https://doi.org/10.1093/imrn/rns124}.

\bibitem[Shinomiya(2016)]{syarticle}
Yoshihiko Shinomiya.
\newblock Veech surfaces and their periodic points.
\newblock \emph{Conformal Geometry and Dynamics of the American Mathematical
  Society}, 20:\penalty0 176--196, 05 2016.
\newblock \doi{10.1090/ecgd/296}.

\bibitem[Hubert and Schmidt(2004)]{Hubert2004INFINITELYGV}
P.~Hubert and Thomas~A. Schmidt.
\newblock Infinitely generated {VEECH} groups.
\newblock \emph{Duke Mathematical Journal}, 123:\penalty0 49--69, 2004.

\bibitem[Kenyon and Smillie(2000)]{Kenyon2000BilliardsOR}
R.~Kenyon and J.~Smillie.
\newblock Billiards on rational-angled triangles.
\newblock \emph{Commentarii Mathematici Helvetici}, 75:\penalty0 65--108, 2000.

\end{thebibliography}





\end{document}